\documentclass[a4paper,11pt]{amsart}

\usepackage{amssymb,cite,url}

\newtheorem{theorem}{Theorem}[section]
\newtheorem{lemma}[theorem]{Lemma}
\newtheorem{corollary}[theorem]{Corollary}
\newtheorem{proposition}[theorem]{Proposition}

\theoremstyle{definition}
\newtheorem{defin}[theorem]{Definition}

\numberwithin{equation}{section}

\makeatletter
\@namedef{subjclassname@2010}{%
\textup{2010} Mathematics Subject Classification}
\makeatother

\DeclareMathOperator{\RE}{Re}
\DeclareMathOperator{\IM}{Im}
\DeclareMathOperator{\Log}{Log}
\DeclareMathOperator{\sgn}{sgn}

\begin{document}

\title[Density function for the Lerch zeta-function]
{The density function for the value-distribution of the Lerch zeta-function and its applications}

\author[M. Mine]{Masahiro Mine}
\address{Department of Mathematics\\ Tokyo Institute of Technology\\ 2-12-1 Ookayama, Meguro-ku, Tokyo 152-8551, Japan}
\email{mine.m.aa@m.titech.ac.jp}

\date{}

\begin{abstract}
The probabilistic study of the value-distributions of zeta-functions is one of the modern topics in analytic number theory. 
In this paper, we study a certain probability measure related to the value-distribution of the Lerch zeta-function. 
We prove that it has a density function which we can explicitly construct. 
Moreover, we prove an asymptotic formula for the number of zeros of the Lerch zeta-function on the right side of the critical line, whose main term is associated with the density function. 
\end{abstract}

%% 11M06	$\zeta (s)$ and $L(s, \chi)$
%% 11M32	Multiple Dirichlet series and zeta functions and multizeta values
%% 11M35	Hurwitz and Lerch zeta functions
%% 11M36	Selberg zeta functions and regularized determinants
%% 11M41	Other Dirichlet series and zeta functions
%% 11R42	Zeta functions and $L$-functions of number fields
%% 11F66	Langlands $L$-functions


\thanks{The work of this work was supported by the Hiki Foundation of Tokyo Institute of Technology and Grant-in-Aid for JSPS Fellows (Grant Number JP19J12037).}

\maketitle

\section{Introduction}\label{s1}

\subsection{Value-distribution of the Riemann zeta-function}\label{s1.1}
Let $\zeta(s)$ be the Riemann zeta-function with $s=\sigma+it$ a complex variable. 
The study of the distribution of values $\zeta(\sigma+it)$ is a classical topic in number theory. 
In 1914, Bohr and Courant \cite{BohrCourant1914} proved that the set $\{ \zeta(\sigma+it) \mid t \in \mathbb{R} \}$ is dense in the complex plane $\mathbb{C}$ for any fixed $1/2<\sigma \leq1$. 
This result was refined by Bohr and Jessen \cite{BohrJessen1930, BohrJessen1932}. 

\begin{theorem}[Bohr--Jessen limit theorem]\label{t1.1}
Let $\sigma>1/2$ be a fixed real number and $R$ be a fixed rectangle in $\mathbb{C}$ whose edges are parallel to the axes. 
Then the limit value 
\begin{equation}\label{e1.1}
W_\sigma(R)
=\lim_{T \to\infty} \frac{1}{T} 
\mu_1\left( \left\{ t \in [0,T] ~\middle|~ \log\zeta(\sigma+it) \in R \right\} \right)
\end{equation}
exists, where $\mu_1$ indicates the one-dimensional Lebesgue measure. 
Moreover, there exists a non-negative real valued continuous function $\mathcal{F}_\sigma(z)$ such that 
\begin{gather*}
W_\sigma(R)
=\int_{R} \mathcal{F}_\sigma(z) \,|dz|,
\end{gather*}
where $|dz|=(2\pi)^{-1}dxdy$ for $z=x+iy$.
\end{theorem}

Here, we determine the blanch of logarithm as follows. 
For $\sigma>1$, we recall that $\zeta(s)$ has the Euler product representation
\begin{equation}\label{e1.2}
\zeta(s)
=\prod_{p} (1-p^{-s})^{-1},
\end{equation}
where $p$ runs through all prime numbers. 
According to this, we define $\log\zeta(s)=-\sum_{p} \Log(1-p^{-s})$ for $\sigma>1$ with $\Log(z)$ the principal blanch. 
In order to extend the definition of $\log\zeta(s)$, we put
\begin{gather*}
G
=\{ \sigma+it \mid \sigma>1/2 \} 
\setminus \bigcup_{\RE(\rho)>1/2} \{ \sigma+i\IM(\rho) \mid 1/2<\sigma \leq \RE(\rho) \}, 
\end{gather*}
where $\rho$ runs through all possible zeros and poles of $\zeta(s)$ with $\RE(\rho)>1/2$. 
Then $\log\zeta(s)$ is defined for $s \in G$ by the analytic continuation along the horizontal path from right. 

Let $\mathcal{B}(\mathbb{C})$ be the class of Borel sets of $\mathbb{C}$. 
For $\sigma>1/2$ and $T>0$, we define a probability measure $P_{\sigma, T}$ on $(\mathbb{C}, \mathcal{B}(\mathbb{C}))$ as 
\begin{gather*}
P_{\sigma, T}(A)
=\frac{1}{T} \mu_1\left( \{ t \in [0,T] \mid \log\zeta(\sigma+it) \in A \} \right).
\end{gather*}
Jessen and Wintner \cite{JessenWintner1935} proved that there exists an absolutely continuous probability measure $P_\sigma$ on $(\mathbb{C}, \mathcal{B}(\mathbb{C}))$ such that $P_{\sigma, T}$ converges weakly to $P_\sigma$ as $T \to\infty$. 
Thus we have 
\begin{equation}\label{e1.3}
\lim_{T \to\infty} P_{\sigma, T}(A)
=P_\sigma(A)
\end{equation}
for any $A \in \mathcal{B}(\mathbb{C})$ whose boundary $\partial A$ satisfies $\mu_2(\partial A)=0$ with the two-dimensional Lebesgue measure $\mu_2$. 
Note that the limit value $W_\sigma(R)$ of Theorem \ref{t1.1} agrees with $P_\sigma(R)$, and the function $\mathcal{F}_\sigma(z)$ is the density of the probability measure $P_\sigma$. 
Borchsenius and Jessen \cite{BorchseniusJessen1948} applied the probabilistic research on values $\log\zeta(\sigma+it)$ to the study of the distribution of $a$-points of $\zeta(s)$. 
Here, we call $z \in \mathbb{C}$ is an $a$-point of $\zeta(s)$ if $z$ satisfies $\zeta(z)=a$. 
If $a \neq 0$, then there exist $a$-points of $\zeta(s)$ in the half plane $\sigma>1/2$ while the Riemann Hypothesis asserts that there are no zeros in this region. 
Borchsenius and Jessen proved that there are a lot of $a$-points of $\zeta(s)$ in the half plane $\sigma>1/2$ when $a \neq0$. 

\begin{theorem}[Borchsenius--Jessen]\label{t1.2}
Let $1/2<\sigma_1<\sigma_2<\infty$ and $a \in \mathbb{C}$. 
Then there exists a limit value 
\begin{equation}\label{e1.4}
C_a(\sigma_1,\sigma_2)
=\lim_{T \to\infty} \frac{1}{T} N_a(T,\sigma_1,\sigma_2),
\end{equation}
where we denote the number of $a$-points of $\zeta(s)$ in the rectangle $\sigma_1<\sigma<\sigma_2$, $0<t<T$ by $N_a(T,\sigma_1,\sigma_2)$. 
Furthermore, the limit value $C_a(\sigma_1,\sigma_2)$ is positive for $\sigma_1 \leq 1$ if $a \neq0$. 
\end{theorem}

Especially, there are infinitely many $a$-points of $\zeta(s)$ in the strip $1/2<\sigma<1$ if $a \neq0$. 
Theorems \ref{t1.1} and \ref{t1.2} have been developed in various ways. 
For example, we consider 
\begin{gather*}
D_{\sigma,T}(R)
=P_{\sigma,T}(R)-P_\sigma(R), 
\end{gather*}
where $\sigma>1/2$, $T>0$, and $R$ is a fixed rectangle as above. 
Matsumoto \cite{Matsumoto1987} obtained an upper bound for $D_{\sigma,T}(R)$, and it was improved by Harman and Matsumoto \cite{HarmanMatsumoto1994}. 
Lamzouri--Lester--Radziwi{\l\l} \cite{LamzouriLesterRadziwill2019} proved the estimate 
\begin{equation}\label{e1.5}
\sup_{R} |D_{\sigma,T}(R)|
=O\left(\frac{1}{(\log{T})^\sigma}\right)
\end{equation}
holds as $T \to\infty$. 
Furthermore, they refined limit formula \eqref{e1.4}. 
In fact, they obtained that the formula 
\begin{equation}\label{e1.6}
N_a(T,\sigma_1,\sigma_2)
=C_a(\sigma_1,\sigma_2)T
+O\left(T \frac{\log\log{T}}{(\log{T})^{\sigma_1/2}}\right)
\end{equation}
holds for $1/2<\sigma_1<\sigma_2<\infty$ and $a \neq0$ as $T \to\infty$. 
Their methods for the proof of \eqref{e1.5} and \eqref{e1.6} were partly based on Guo \cite{Guo1996a, Guo1996b}, who studied the value-distribution of $(\zeta'/\zeta)(s)$ instead of $\log\zeta(s)$. 

In \cite{Matsumoto1989, Matsumoto1990, Matsumoto1992}, Matsumoto generalized Theorem \ref{t1.1} in a class of $L$-functions possessing polynomial Euler products. 
Then Ihara \cite{Ihara2008} named density functions such as $\mathcal{F}_\sigma(z)$ ``$M$-functions'' for the value-distributions of $L$-functions. 
Ihara and Matsumoto \cite{IharaMatsumoto2011a, IharaMatsumoto2011b, IharaMatsumoto2014} studied $M$-functions for various value-distributions of Dirichlet $L$-functions, Hecke $L$-functions, and so on; see the survey of Matsumoto \cite{Matsumoto2019}. 
In this paper, we study the $M$-function for the value-distribution of the Lerch zeta-function $L(\lambda,\alpha,s)$, which does not have the Euler product in general.

\subsection{Hurwitz and Lerch zeta-functions}\label{s1.2}
Let $\lambda \in \mathbb{R}$ and $0<\alpha \leq1$. 
Then we define the Hurwitz zeta-function $\zeta(s, \alpha)$ as 
\begin{gather*}
\zeta(s, \alpha)
=\sum_{n=0}^{\infty} \frac{1}{(n+\alpha)^s}
\end{gather*}
for $\sigma>1$, which is continued holomorphically to the whole complex plane except for a simple pole at $s=1$. 
Since we have $\zeta(s,1)=\zeta(s)$ by definition, the Hurwitz zeta-function is a generalization of the Riemann zeta-function. 
The Lerch zeta-function 
\begin{gather*}
L(\lambda,\alpha,s)
=\sum_{n=0}^{\infty} \frac{e^{2\pi i \lambda n}}{(n+\alpha)^s}
\end{gather*}
is a further generalization of $\zeta(s)$. 
Indeed, it deduces the Hurwitz zeta-function if $\lambda \in \mathbb{Z}$. 
Some of the results on $\zeta(s)$ are naturally generalized to $\zeta(s,\alpha)$ and $L(\lambda,\alpha,s)$. 
For instance, let $N(T;\lambda,\alpha)$ be the number of non-trivial zeros $\rho$ of $L(\lambda, \alpha, s)$ with $0<\IM(\rho)<T$. 
Then we have
\begin{gather*}
N(T;\lambda,\alpha)
=\frac{T}{2\pi} \log\frac{T}{2\pi e \alpha \lambda}+O(\log{T})
\end{gather*}
for every $0<\lambda \leq1$ and $0<\alpha \leq1$; see \cite[Theorem 6]{GarunkstisLaurincikas1999}. 
It is a generalization of the classical result of von-Mangoldt on $N(T)=N(T;1,1)$. 
On the other hand, a significant difference between $\zeta(s)$ and $\zeta(s,\alpha)$ or $L(\lambda,\alpha,s)$ arises when we consider their Euler products and zeros off the critical line $\sigma=1/2$. 
By definition of $L(\lambda,\alpha,s)$, we have 
\begin{gather*}
L(1,1/2,s)
=(2^s-1) \zeta(s) 
\quad\text{and}\quad 
L(1/2,1,s)
=(1-2^{1-s}) \zeta(s).
\end{gather*}
Hence $L(1,1,s)$, $L(1,1/2,s)$, and $L(1/2,1,s)$ have no zeros in the half plane $\sigma>1$ due to the Euler product of the Riemann zeta-function \eqref{e1.2}. 
However, $\zeta(s,\alpha)$ and $L(\lambda,\alpha,s)$ do not have the Euler product representation in general. 
For example, we see that 
\begin{gather*}
\zeta(s,1/3)
=3^s \sum_{n=1}^{\infty} \frac{\delta_{\equiv1(\bmod 3)}(n)}{n^s},
\end{gather*}
where $\delta_{\equiv1(\bmod 3)}(n)$ is equal to $1$ if $n \equiv1\,(\bmod\,3)$ and $0$ otherwise, which is not multiplicative. 
Because of the lack of Euler products, one may hope that Hurwitz and Lerch zeta-functions have zeros in the half plane $\sigma>1$. 
Davenport and Heilbronn \cite{DavenportHeilbronn1936a} proved that there are infinitely many zeros of the Hurwitz zeta-function $\zeta(s,\alpha)$ in $\sigma>1$ if $\alpha \neq1,1/2$ is a rational or transcendental number. 
Moreover, Cassels \cite{Cassels1961} proved the same result in the remaining case, i.e. $\zeta(s,\alpha)$ has infinitely many zeros in $\sigma>1$ even if $\alpha$ is algebraic irrational. 
We then proceed to considering zeros in the strip $1/2<\sigma \leq1$. 
Although it is conjectured that $\zeta(s)$ has no zeros there, $\zeta(s,\alpha)$ has infinitely many zeros in $1/2<\sigma \leq1$ if $\alpha \neq1,1/2$ is either rational or transcendental. 
This result is a consequence of a property of $\zeta(s,\alpha)$ so-called the strong universality proved by Bagchi \cite{Bagchi1981} and Gonek \cite{Gonek1979}. 
Unfortunately, such a property has not yet been known in the case that $\alpha$ is an algebraic irrational number. 

Garunk\v{s}tis and Laurin\v{c}ikas proved similar results on zeros of the Lerch zeta-function on the right side of the critical line; see \cite[Chapter 8]{LaurincikasGarunkstis2002}. 
Let $N(T, \sigma_1, \sigma_2; \lambda, \alpha)$ be the number of zeros of $L(\lambda,\alpha,s)$ in the rectangle $\sigma_1<\sigma<\sigma_2$, $0<t<T$. 
Then, combining the results of Garunk\v{s}tis and Laurin\v{c}ikas, we obtain the following result. 

\begin{theorem}[Garunk\v{s}tis--Laurin\v{c}ikas]\label{t1.3}
Let $\lambda \in \mathbb{R}$ and $0<\alpha \leq1$. 
Assume that $\alpha$ is a transcendental number. 
Then there exists a constant $\delta(\alpha)>0$ depending only on $\alpha$ such that, for any fixed real numbers $\sigma_1$, $\sigma_2$ with $1/2<\sigma_1<\sigma_2<1+\delta(\alpha)$, there exist positive constants $C_1=C_1(\sigma_1, \sigma_2; \lambda, \alpha)$ and $C_2=C_2(\sigma_1, \sigma_2; \lambda, \alpha)$ such that 
\begin{gather*}
C_1 T 
<N(T, \sigma_1, \sigma_2; \lambda, \alpha)
<C_2 T
\end{gather*}
holds for sufficiently large $T$.
\end{theorem}

In this paper, we study the value-distribution of $L(\lambda,\alpha,s)$ and the asymptotic behavior of $N(T, \sigma_1, \sigma_2; \lambda, \alpha)$ according to analogues of Theorems \ref{t1.1} and \ref{t1.2} for $L(\lambda,\alpha,s)$. 
Let $P_{\sigma, T}(~\cdot~; \lambda, \alpha)$ be a probability measure on $(\mathbb{C}, \mathcal{B}(\mathbb{C}))$ defined as 
\begin{gather*}
P_{\sigma, T}(A; \lambda, \alpha)
=\frac{1}{T} \mu_1\left( \{ t \in [0,T] \mid L(\lambda, \alpha, \sigma+it) \in A \} \right).
\end{gather*}
Then Garunk\v{s}tis and Laurin\v{c}ikas \cite{GarunkstisLaurincikas1996} proved the following result. 

\begin{theorem}[Garunk\v{s}tis--Laurin\v{c}ikas]\label{t1.4}
Let $\sigma>1/2$, $\lambda \in \mathbb{R}$, and $0<\alpha \leq1$. 
Then there exists a probability measure $P_\sigma(~\cdot~; \lambda, \alpha)$ on $(\mathbb{C}, \mathcal{B}(\mathbb{C}))$ such that $P_{\sigma, T}(~\cdot~; \lambda, \alpha)$ converges weakly to $P_\sigma(~\cdot~; \lambda, \alpha)$ as $T \to\infty$. 
\end{theorem}

Note that Theorem \ref{t1.4} holds with an arbitrary $0<\alpha \leq1$, but it does not ensure the absolutely continuity of the limit measure $P_\sigma(~\cdot~; \lambda, \alpha)$. 
If $P_\sigma(~\cdot~; \lambda, \alpha)$ is absolutely continuous, we obtain
\begin{equation}\label{e1.7}
\lim_{T \to\infty} P_{\sigma, T}(A; \lambda, \alpha)
=P_\sigma(A; \lambda, \alpha)
\end{equation}
for any $A \in \mathcal{B}(\mathbb{C})$ with $\mu_2(\partial A)=0$ as we have seen in \eqref{e1.3}. 
However, Theorem \ref{t1.3} tells us only that \eqref{e1.7} is valid for $A \in \mathcal{B}(\mathbb{C})$ with $P_\sigma(\partial A;\lambda,\alpha)=0$ without the absolutely continuity. 

We have assumed that the parameter $\alpha$ is transcendental in Theorem \ref{t1.3}. 
In the study of the value-distributions of $\zeta(s,\alpha)$ and $L(\lambda,\alpha,s)$, this assumption sometimes can be replaced with the assumption that the system $\{\log(n+\alpha)\}_{n \geq0}$ is linearly independent over $\mathbb{Q}$. 
Then, we obtain the following result as an analogue of Theorem \ref{t1.2} under the latter assumption. 

\begin{theorem}\label{t1.5}
Let $\lambda \in \mathbb{R}$ and $0<\alpha \leq1$. 
Assume that the system $\{\log(n+\alpha)\}_{n \geq0}$ is linearly independent over $\mathbb{Q}$. 
Then the limit value 
\begin{equation}\label{e1.8}
C(\sigma_1, \sigma_2;\alpha)
=\lim_{T \to\infty} \frac{1}{T} N(T, \sigma_1, \sigma_2; \lambda, \alpha)
\end{equation}
exists for any $1/2<\sigma_1<\sigma_2<\infty$. 
\end{theorem}

It is notable that the limit value is determined only from $\sigma_1$, $\sigma_2$, $\alpha$ and does not depend on $\lambda$ although the constants $C_1$ and $C_2$ of Theorem \ref{t1.3} may depend on $\lambda$. 
One can prove Theorem \ref{t1.5} by just adapting the method of Borchsenius--Jessen \cite{BorchseniusJessen1948}. 
However, we give a full proof in Appendix of this paper since the author could not find a suitable reference for the proof.

\subsection{Statements of results}\label{s1.3}
Comparing the results on $\log\zeta(s)$ in Section \ref{s1.1} and on $L(\lambda,\alpha,s)$ in Section \ref{s1.2}, we hope that there is room for improvement of limit formulas \eqref{e1.7} and \eqref{e1.8} as well as \eqref{e1.5} and \eqref{e1.6}. 
The main purpose of this paper is to refine Theorems \ref{t1.4} and \ref{t1.5} in this sense. 
For this, we need a further restriction to the parameter $\alpha$. 

\begin{defin}\label{d1.6}
We define the \emph{admissible subset} $\mathfrak{S}$ as the collection of all $\alpha \in (0,1]$ which satisfy the following conditions \eqref{c1} and \eqref{c2}. 
\begin{enumerate}
\item 
\emph{Independence.} \label{c1}
The system $\{\log(n+\alpha)\}_{n \geq0}$ is linearly independent over $\mathbb{Q}$. 
\item 
\emph{Spacing.} \label{c2}
There exists a constant $\Omega(\alpha)>0$ such that for large $X>0$ and for any positive integer $N$, we have 
\begin{gather*}
\left|\sum_{j=1}^{N} \epsilon_j \log(n_j+\alpha)\right|
\geq X^{-\Omega(\alpha)N^2}
\end{gather*}
provided $\sum_{j=1}^{N} \epsilon_j \log(n_j+\alpha) \neq0$ with $0 \leq n_1, \ldots, n_N \leq X$ arbitrary integers and $\epsilon_1, \ldots, \epsilon_N \in \{\pm1\}$. 
\end{enumerate}
\end{defin}

We see that almost all $\alpha \in (0,1]$ belong to $\mathfrak{S}$ with respect to $\mu_1$ by applying the results of transcendental number theory. 
For the proof of this fact and more information about the admissible subset $\mathfrak{S}$; see Section \ref{s2.1}. 
Then, we prove the following result as an improvement of Theorem \ref{t1.4}. 

\begin{theorem}\label{t1.7}
Let $\lambda \in \mathbb{R}$ and $0<\alpha \leq1$. 
Let $\sigma_1>0$ be a large fixed real number. 
Let $\epsilon_1>0$ be a small fixed real number. 
Assume that $\alpha$ is a member of $\mathfrak{S}$. 
Then for each $\epsilon>0$, there exists a positive constant $T_0=T_0(\alpha, \sigma_1, \epsilon_1, \epsilon)$ such that for all $1/2+\epsilon_1 \leq \sigma \leq \sigma_1$ and for all $T \geq T_0$, we have
\begin{gather*}
P_{\sigma, T}(R; \lambda, \alpha)
=P_\sigma(R; \lambda, \alpha)
+O\left( (\mu_2(R)+1) (\log{T})^{-1/4+\epsilon} \right), 
\end{gather*}
where $R$ is any rectangle in $\mathbb{C}$ whose edges are parallel to the axes and have length greater than $(\log{T})^{-1/4+\epsilon}$. 
The implied constant depends only on $\alpha$, $\lambda$, $\sigma_1$, and $\epsilon_1$. 
Moreover, there exists a non-negative smooth function $M_\sigma(z; \alpha)$ on $\mathbb{C}$ such that we have 
\begin{gather*}
P_\sigma(A; \lambda, \alpha)
=\int_{A} M_\sigma(z; \alpha) \,|dz|
\end{gather*}
for every $\sigma>1/2$ and $A \in \mathcal{B}(\mathbb{C})$. 
\end{theorem}

The second statement of Theorem \ref{t1.7} assert that the limit measure $P_\sigma(~\cdot~; \lambda, \alpha)$ is absolute continuous with the density $M_\sigma(z; \alpha)$ if $\alpha \in \mathfrak{S}$. 
The function $M_\sigma(z; \alpha)$ is an analogue of $\mathcal{F}_\sigma(z)$ of Theorem \ref{t1.1}. 
We again remark that it does not depend on $\lambda$. 
The final result of this paper is as follows, which is an improvement of Theorems \ref{t1.3} and \ref{t1.5} in the case that $\alpha$ belongs to the admissible subset $\mathfrak{S}$. 

\begin{theorem}\label{t1.8}
Let $1/2<\sigma_1<\sigma_2<\infty$, $\lambda \in \mathbb{R}$, and $\alpha \in \mathfrak{S}$. 
Then there exists an absolute constant $A>0$ such that  
\begin{gather*}
N(T, \sigma_1, \sigma_2; \lambda, \alpha)
=C(\sigma_1, \sigma_2; \alpha)T+O\left(T(\log{T})^{-A}\right)
\end{gather*}
as $T \to\infty$, where the implied constant depends only on $\sigma_1$, $\sigma_2$, $\lambda$ and $\alpha$. 
Furthermore, the constant $C(\sigma_1, \sigma_2; \alpha)$ is given by 
\begin{gather*}
C(\sigma_1, \sigma_2; \alpha)
=\frac{1}{2\pi} \int_{\sigma_1}^{\sigma_2} 
\left( \int_{\mathbb{C}} \log |z| \frac{\partial^2}{\partial \sigma^2} M_\sigma(z; \alpha) \,|dz| \right) \,d\sigma. 
\end{gather*}
\end{theorem}

This result is also an analogue of the result of Guo \cite[Theorem 1.1.2]{Guo1996b} on the distribution of zeros of $\zeta'(s)$. 
We adapt the method of Guo for a part of the proof of Theorem \ref{t1.8}. 

\vspace{\baselineskip}

This paper consists of seven sections. 
Section \ref{s2} contains preliminaries for the proofs of the main results. 
In Section \ref{s2.1}, we study the admissible subset $\mathfrak{S}$ more precisely. 
We recall the notion of $S$-numbers introduced by Mahler \cite{Mahler1932b} and show that every $S$-number in $(0,1]$ belongs to $\mathfrak{S}$. 
In Section \ref{s2.2}, we construct the density function $M_\sigma(z; \alpha)$. 
For this, we review the theory of equidistributions on circles studied by Jessen and Wintner \cite{JessenWintner1935}. 
The function $M_\sigma(z; \alpha)$ is defined as an infinite convolution of such distributions. 
In Section \ref{s3}, we consider several mean values of the Lerch zeta-function. 
We prove Theorem \ref{t3.2} in this section which assert that the mean values are expressed as the integrals involving the density function $M_\sigma(z; \alpha)$. 
Both of the proofs of Theorems \ref{t1.7} and \ref{t1.8} are based on Theorem \ref{t3.2}. 
Theorem \ref{t1.7} is proved in Section \ref{s4}, where we use the Beurding--Selberg functions which approximate well the signum function $\sgn(x)$ on $\mathbb{R}$. 
In order to prove Theorem \ref{t1.8}, we need to examine the analytic properties of the function $M_\sigma(z; \alpha)$ more precisely. 
Further results on $M_\sigma(z; \alpha)$ are given in Section \ref{s5}. 
Finally, we prove Theorem \ref{t1.8} in Section \ref{s6}. 
A mean value theorem for $\log|L(\lambda, \alpha, s)|$ is a key for the proof. 
The last section is an appendix. 
We give a proof of Theorem \ref{t1.5} in this section by following the method of Borchsenius--Jessen \cite{BorchseniusJessen1948}. 

\vspace{\baselineskip}

Throughout this paper, we identify the complex plane $\mathbb{C}$ with $\mathbb{R}^2$ by the bijection $x+iy \mapsto (x,y)$. 
Hence we regard functions on $\mathbb{C}$ as functions on $\mathbb{R}^2$, and for example, we write the two-dimensional $L^p$-space by $L^p(\mathbb{C})=L^p(\mathbb{R}^2)$. 
We also use the symbols $\ll$ and $\gg$ as the usual Vinogradov symbols.

\section{Preliminaries}\label{s2}

\subsection{Admissible subset $\mathfrak{S}$}\label{s2.1}
Definition \ref{d1.6} is motivated by transcendental number theory. 
In fact, we find that any transcendental number satisfies condition \eqref{c1}. 
For condition \eqref{c2}, we recall the notion of $S$-numbers introduced by Mahler \cite{Mahler1932b}. 
According to the notation of the book of Baker \cite{Baker1990}, for a complex number $\xi$, and for positive integers $n$ and $h$, let $P(x) \in \mathbb{Z}[x]$ be a polynomial with degree at most $n$ and height at most $h$ for which $|P(\xi)|$ takes the smallest non-zero value. 
Here the height of a polynomial is the maximal value of the absolute values of the coefficients. 
Then, we determine a real number $\omega_{n, h}(\xi) \in (0,\infty)$ by the equation $|P(\xi)|=h^{-n \omega_{n, h}(\xi)}$. 
We further define $\omega(\xi) \in [0,\infty]$ as
\begin{gather*}
\omega(\xi)
=\limsup_{n \to\infty} \limsup_{h \to\infty} \omega_{n, h}(\xi).
\end{gather*}

\begin{defin}\label{d2.1}
A complex number $\xi$ is called an \emph{$S$-number} if $0<\omega(\xi)<\infty$. 
\end{defin}

Hence, if $\xi$ is an $S$-number, there exists a positive real number $\omega'(\xi)$ such that we have 
\begin{equation}\label{e2.1}
|P(\xi)|
\geq h^{-n \omega'(\xi)}
\end{equation}
for every polynomial $P(x) \in \mathbb{Z}[x]$ with degree at most $n$ and height at most $h$ if $P(\xi) \neq0$. 
Then, the following result holds. 

\begin{lemma}\label{l2.2}
Every $S$-number in $(0,1]$ belongs to the admissible subset $\mathfrak{S}$. 
\end{lemma}

\begin{proof}
It is known that all $S$-numbers are transcendental \cite[Section 8.2]{Baker1990}, and therefore every $S$-number $\alpha \in (0,1]$ satisfies condition \eqref{c1}. 
Then we prove that it also satisfies condition \eqref{c2}. 
Let $\mu$ and $\nu$ be non-negative integers with $\mu+\nu=N$, where $N \geq1$. 
Then we take $N$ integers $0 \leq m_1,\ldots, m_\mu, n_1,\dots,n_\nu \leq X$ which satisfy $(m_1+\alpha) \cdots (m_\mu+\alpha) \neq (n_1+\alpha) \cdots (n_\nu+\alpha)$. 
Here, if one of $\mu$ and $\nu$ is equal to $0$, the product $(m_1+\alpha) \cdots (m_\mu+\alpha)$ or $(n_1+\alpha)\cdots (n_\nu+\alpha)$ is interpreted as $1$. 
We have 
\begin{align}\label{e2.2}
&\left|\sum_{j=1}^{\mu} \log(m_j+\alpha) - \sum_{k=1}^{\nu} \log(n_k+\alpha)\right|
=\left|\log \frac{(m_1+\alpha) \cdots (m_\mu+\alpha)}{(n_1+\alpha) \cdots (n_\nu+\alpha)}\right| \\
&\geq \frac{|(m_1+\alpha) \cdots (m_\mu+\alpha)-(n_1+\alpha) \cdots (n_\nu+\alpha)|}
{\max\{ (m_1+\alpha) \cdots (m_\mu+\alpha), (n_1+\alpha) \cdots (n_\nu+\alpha) \}} \nonumber \\
&\geq \frac{|P(\alpha)|}{(2X)^N} \nonumber
\end{align}
with the polynomial
\begin{gather*}
P(x)
=(m_1+x) \cdots (m_\mu+x)-(n_1+x) \cdots (n_\nu+x).
\end{gather*}
The degree of $P(x)$ is at most $\mu+\nu=N$, and the height is at most 
\begin{gather*}
\binom{N}{\lfloor N/2\rfloor} X^N
\leq (2X)^N, 
\end{gather*}
where $\lfloor x \rfloor$ is the maximum integer not greater than $x$. 
Hence, by \eqref{e2.1}, we have 
\begin{gather*}
|P(\alpha)|
\geq (2X)^{-\omega'(\alpha) N^2}.
\end{gather*}
Therefore we obtain the desired inequality by \eqref{e2.2} when $X$ is large. 
\end{proof}

Almost all real numbers are $S$-numbers with respect to $\mu_1$; see \cite[Theorem 8.2]{Baker1990}. 
Hence we deduce the following corollary. 

\begin{corollary}\label{c2.3}
Almost all $\alpha$ in $(0,1]$ belong to $\mathfrak{S}$ with respect to $\mu_1$. 
\end{corollary}

A typical example of $S$-numbers is Napier's constant $e$, and we see that its fractional part $\{ e \}$ is a member of $\mathfrak{S}$. 
On the other hand, Liouville numbers such as $\sum_{k} 10^{-k!}$ are transcendental but not $S$-numbers. 
However, it seems difficult to see Liouville numbers in $(0,1]$ are belongs to $\mathfrak{S}$ or not.

\subsection{Equidistributions on circles}\label{s2.2}
According to \cite[Section 5]{JessenWintner1935}, let $S$ be the circle $|z|=r$ in the complex plane with $r>0$, and define the equidistribution on $S$ as the probability measure $\phi$ on $(\mathbb{C}, \mathcal{B}(\mathbb{C}))$ given by
\begin{gather*}
\phi(A)
=\widetilde{\mu}_S(A \cap S),
\end{gather*}
where $\widetilde{\mu}_S$ is the normalized Haar measure of $S$. 
For $0<\alpha \leq1$, let $S_n$ be the circles $|z|=(n+\alpha)^{-\sigma}$ for $n \geq0$. 
Denote the equidistributions on $S_n$ by $\phi_n$. 
Then we consider the infinite convolution $\phi=\phi_0*\phi_1*\cdots$. 

\begin{proposition}\label{p2.4}
For any $\sigma>1/2$, the convolution measure $\phi_0*\phi_1*\cdots*\phi_n$ converses weakly to a probability measure $\phi$ as $n \to\infty$. 
Moreover, the limit measure $\phi$ is absolutely continuous with a continuous density function $M_\sigma(z;\alpha)$, i.e. there exists a non-negative real valued continuous function $M_\sigma(z;\alpha)$ such that 
\begin{gather*}
\phi(A)=\int_A M_\sigma(z;\alpha)\,|dz|, 
\qquad 
A \in \mathcal{B}(\mathbb{C}). 
\end{gather*}
\end{proposition}

\begin{proof}
This is a simple application of \cite[Theorem 10]{JessenWintner1935}. 
\end{proof}

In general, we define the Fourier transform
\begin{gather*}
\widetilde{f}(z)
=\int_{\mathbb{C}} f(w) \psi_z(w) \,|dw|
\end{gather*}
for $f \in L^1(\mathbb{C})$ with $\psi_z(w)=\exp\left(i \RE(z \overline{w})\right)$. 
Then we have 
\begin{equation}\label{e2.3}
\widetilde{M}_\sigma(z; \alpha)
=\int_{\mathbb{C}} M_\sigma(w; \alpha) \psi_z(w) \,|dw|
=\int_{\mathbb{C}} \psi_z(w) \,d\phi(w).
\end{equation}
Since $\phi=\phi_0*\phi_1*\cdots$, we find that 
\begin{gather*}
\int_{\mathbb{C}} \psi_z(w) \,d\phi(w)
=\prod_{n=0}^{\infty} \int_{\mathbb{C}} \psi_z(w) \,d\phi_n(w)
\end{gather*}
by the argument in \cite[Section 4]{JessenWintner1935}. 
Furthermore, by \cite[Section 5]{JessenWintner1935}, we have 
\begin{gather*}
\int_{\mathbb{C}} \psi_z(w) \,d\phi_n(w)
=J_0(|z|(n+\alpha)^{-\sigma}),
\end{gather*}
where $J_0(x)$ is the Bessel function of the first kind of order $0$. 
Therefore, the Fourier transform $\widetilde{M}_\sigma(z; \alpha)$ is expressed as the infinite product
\begin{equation}\label{e2.4}
\widetilde{M}_\sigma(z; \alpha)
=\prod_{n=0}^{\infty} J_0(|z|(n+\alpha)^{-\sigma}). 
\end{equation}
More detailed properties on the function $M_\sigma(z; \alpha)$ are as follows. 

\begin{proposition}\label{p2.5}
We have the following properties on $M_\sigma(z; \alpha)$. 
\begin{itemize}
\item[$(\mathrm{i})$]
For any $\sigma>1/2$ and $0<\alpha \leq1$, the function $M_\sigma(z; \alpha)$ possesses continuous partial derivative of arbitrarily high order as a function of $x$ and $y$ with $z=x+iy$. 
Their growths are estimated as for any $c>0$ and for $\sigma \geq 1/2+\epsilon_1$ with small $\epsilon_1>0$, 
\begin{equation}\label{e2.5}
\frac{\partial^{m+n}}{\partial x^m \partial y^n} M_\sigma(z; \alpha)
\ll e^{-c|z|^2},
\qquad 
|z| \to\infty 
\end{equation}
for all non-negative integers $m$ and $n$, where the implied constant depends only on $\alpha$, $m$, $n$, and $\epsilon_1$. 
\item[$(\mathrm{ii})$]
For any $\sigma>1/2$ and $0<\alpha \leq1$, we have 
\begin{equation}\label{e2.6}
\int_{\mathbb{C}} M_\sigma(z; \alpha) \,|dz|
=1.
\end{equation}
\item[$(\mathrm{iii})$]
If $1/2<\sigma \leq1$, then $M_\sigma(z; \alpha)>0$ for all $z \in \mathbb{C}$. 
\item[$(\mathrm{iv})$]
$\widetilde{M}_\sigma(z; \alpha)$ is a real valued function with $|\widetilde{M}_\sigma(z; \alpha)| \leq1$. 
The values $\widetilde{M}_\sigma(z; \alpha)$ are determined only from $\sigma$, $\alpha$, and $|z|$, i.e. $\widetilde{M}_\sigma(z_1; \alpha)=\widetilde{M}_\sigma(z_2; \alpha)$ holds if $|z_1|=|z_2|$.
\end{itemize}
\end{proposition}

\begin{proof}
\begin{itemize}
\item[$(\mathrm{i})$]
The former statement is deduced from \cite[Theorem 10]{JessenWintner1935}, and the latter \eqref{e2.5} is due to \cite[Theorem 16]{JessenWintner1935}. 
\item[$(\mathrm{ii})$]
By \eqref{e2.3} and \eqref{e2.4}, we have
\begin{gather*}
\int_{\mathbb{C}} M_\sigma(z; \alpha) \,|dz|
=\widetilde{M}_\sigma(0; \alpha)
=\prod_{n=0}^{\infty} J_0(0).
\end{gather*}
Thus we obtain \eqref{e2.6} since $J_0(0)=1$. 
\item[$(\mathrm{iii})$]
This is again deduced from \cite[Theorem 10]{JessenWintner1935}. 
\item[$(\mathrm{iv})$]
We know that $J_0(x)$ is real valued for $x \in \mathbb{R}$ and satisfies $|J_0(x)| \leq1$. 
Hence formula \eqref{e2.4} yields the first statement. 
The remaining statement is also deduced from \eqref{e2.4}.
\end{itemize}
\end{proof}

\section{Mean values of the Lerch zeta-function}\label{s3}
In this section, we study the mean values  
\begin{equation}\label{e3.1}
\frac{1}{T} \int_{0}^{T} \Phi(L(\lambda, \alpha, \sigma+it)) \,dt
\end{equation}
for test functions $\Phi(z)$ via the density function $M_\sigma(z; \alpha)$ defined in Section \ref{s2.2}. 

First, we define two classes of test functions. 
Let $\mathcal{S}=\mathcal{S}(\mathbb{C})=\mathcal{S}(\mathbb{R}^2)$ denote the Schwartz space on $\mathbb{C}$, which consists of all rapidly decreasing $C^\infty$-functions whose partial derivatives of arbitrarily high order also rapidly decrease. 
Next, according to \cite[Section 9]{IharaMatsumoto2011a}, we also define the class $\Lambda$ as the set of all functions $f \in L^1(\mathbb{C}) \cap L^\infty(\mathbb{C})$ with $\widetilde{f} \in L^1(\mathbb{C}) \cap L^\infty(\mathbb{C})$ that satisfy the inverse formula 
\begin{gather*}
f(z)
=\int_{\mathbb{C}} \widetilde{f}(w) \psi_{-z}(w) \,|dw|.
\end{gather*}
Note that the Schwartz space $\mathcal{S}$ is included in the class $\Lambda$, and especially, any compactly supported $C^\infty$-function belongs to the class $\Lambda$. 
By Proposition \ref{p2.5}, the function $M_\sigma(z; \alpha)$ belongs to $\mathcal{S}$, and hence its Fourier transform $\widetilde{M}_\sigma(z; \alpha)$ does. 

Then, we state the following two results on the mean values \eqref{e3.1}. 

\begin{theorem}\label{t3.1}
Let $\lambda \in \mathbb{R}$ and $\alpha \in \mathfrak{S}$. 
Let $\sigma_1$ be a large fixed positive constant, and let $\theta, \delta>0$ with $\theta+\delta<1/4$. 
Then there exists $T_0=T_0(\alpha, \sigma_1, \theta, \delta)>0$ such that for all $T \geq T_0$ and for all $\sigma \in [1/2+(\log{T})^{-\theta}, \sigma_1]$, we have 
\begin{gather*}
\frac{1}{T} \int_{0}^{T} \psi_z\left(L(\lambda, \alpha, \sigma+it)\right) \,dt
=\widetilde{M}_\sigma(z; \alpha)
+O\left(\exp\left(-\frac{1}{2}(\log{T})^{\theta/2}\right)\right)
\end{gather*}
for any $z \in \Omega$, where 
\begin{gather*}
\Omega
=\{ x+iy \in \mathbb{C} \mid |x| \leq (\log{T})^\delta, |y| \leq (\log{T})^\delta \}.
\end{gather*}
The implied constant depends only on $\lambda$ and $\sigma_1$. 
\end{theorem}

\begin{theorem}\label{t3.2}
Let $\lambda \in \mathbb{R}$ and $\alpha \in \mathfrak{S}$. 
Let $\sigma_1$ be a large fixed positive constant, and let $\theta, \delta>0$ with $\theta+\delta<1/4$. 
Then there exists $T_0=T_0(\alpha, \sigma_1, \theta, \delta)>0$ such that for all $T \geq T_0$ and for all $\sigma \in [1/2+(\log{T})^{-\theta}, \sigma_1]$, we have 
\begin{gather*}
\frac{1}{T} \int_{0}^{T} \Phi(L(\lambda, \alpha, \sigma+it)) \,dt
=\int_{\mathbb{C}} \Phi(z) M_\sigma(z; \alpha) \,|dz| 
+E
\end{gather*}
for any function $\Phi$ in the class $\Lambda$, where $E$ is estimated as 
\begin{gather*}
E
\ll \exp\left(-\frac{1}{2}(\log{T})^{\theta/2}\right) \int_{\Omega} \left|\widetilde{\Phi}(z)\right| \,|dz|
+\int_{\mathbb{C} \setminus \Omega} \left|\widetilde{\Phi}(z)\right| \,|dz|.
\end{gather*}
The implied constant depends only on $\lambda$ and $\sigma_1$. 
\end{theorem}

Theorem \ref{t3.2} is an analogue of \cite[Theorem 1.1.1]{Guo1996a} for Lerch zeta-functions. 
We first prove that Theorem \ref{t3.1} implies Theorem \ref{t3.2}. 

\begin{proof}[Proof of Theorem \ref{t3.2} assuming Theorem \ref{t3.1}]
By definition of the class $\Lambda$, we have 
\begin{gather*}
\Phi(L(\lambda, \alpha, s))
=\int_{\mathbb{C}} \widetilde{\Phi}(z) \psi_{-z}(L(\lambda, \alpha, s)) \,|dz|
\end{gather*}
for any $\Phi \in \Lambda$.
Then we use Theorem \ref{t3.1}. 
For all $T \geq T_0$ and $\sigma \in [1/2+(\log{T})^{-\theta}, \sigma_1]$, we have
\begin{align*}
\frac{1}{T} \int_{0}^{T} \Phi(L(\lambda, \alpha, \sigma+it)) \,dt
&=\int_{\mathbb{C}} \widetilde{\Phi}(z) 
\frac{1}{T} \int_{0}^{T} \psi_{-z}(L(\lambda, \alpha, \sigma+it)) \,dt |dz|\\
&=\int_{\Omega} \widetilde{\Phi}(z) 
\frac{1}{T} \int_{0}^{T} \psi_{-z}(L(\lambda, \alpha, \sigma+it)) \,dt |dz| 
+E_1\\
&=\int_{\Omega} \widetilde{\Phi}(z) \widetilde{M}_\sigma(-z; \alpha) \,|dz| 
+E_1+E_2 \\
&=\int_{\mathbb{C}} \widetilde{\Phi}(z) \widetilde{M}_\sigma(-z; \alpha) \,|dz|
+E_1+E_2+E_3.
\end{align*}
The above error terms are estimated as 
\begin{gather*}
E_1,~E_3 
\ll \int_{\mathbb{C} \setminus \Omega} \left|\widetilde{\Phi}(z)\right| \,|dz| \\
E_2
\ll_{\lambda,\sigma_1} \exp\left(-\frac{1}{2}(\log{T})^{\theta/2}\right) \int_{\Omega} \left|\widetilde{\Phi}(z)\right| \,|dz|
\end{gather*}
due to the inequalities $|\psi_z(w)| \leq1$ and $|\widetilde{M}_\sigma(z; \alpha)| \leq1$.
Finally, we have 
\begin{gather*}
\int_{\mathbb{C}} \widetilde{\Phi}(z) \widetilde{M}_\sigma(-z; \alpha) \,|dz|
=\int_{\mathbb{C}} \widetilde{\Phi}(z) \widetilde{M}_\sigma(z; \alpha) \,|dz|
=\int_{\mathbb{C}} \Phi(z) M_\sigma(z; \alpha) \,|dz| 
\end{gather*}
by $\mathrm{(iv)}$ of Proposition \ref{p2.5} and Parseval's identity. 
Hence the result follows. 
\end{proof}

In the remaining part of this section, we prove Theorem \ref{t3.1}. 
The proof consists of the proofs of the following four Propositions \ref{p3.3}, \ref{p3.5}, \ref{p3.6}, and \ref{p3.7}. 

\begin{proposition}\label{p3.3}
Let $\lambda \in \mathbb{R}$ and $0<\alpha \leq1$. 
Let $\sigma_1$ be a large fixed positive constant, and let $\theta, \delta>0$ with $\theta<2/3$. 
Then there exists $T_0=T_0(\theta)>0$ such that for all $T \geq T_0$ and for all $\sigma \in [1/2+(\log{T})^{-\theta},  \sigma_1]$, we have 
\begin{gather*}
\frac{1}{T} \int_{0}^{T} \psi_z(L(\lambda, \alpha, \sigma+it)) \,dt
=\frac{1}{T} \int_{0}^{T} \psi_z\left(\sum_{0 \leq n \leq X} \frac{e^{2\pi i \lambda n}}{(n+\alpha)^{\sigma+it}}\right) \,dt
+E_1
\end{gather*}
for any $z \in \Omega$, where $X=\exp((\log{T})^{3\theta/2})$, and $E_1$ is estimated as 
\begin{gather*}
E_1
\ll_{\lambda,\sigma_1} \exp\left(-\frac{1}{2}(\log{T})^{\theta/2}\right).
\end{gather*}
\end{proposition}

To prove Proposition \ref{p3.3}, we use the following lemma. 

\begin{lemma}\label{l3.4}
Let $0<\lambda \leq1$ and $0<\alpha \leq1$. 
Then, for any $\sigma>1/2$ and $2\pi \leq |t| \leq \pi \lambda Y$, we have 
\begin{gather*}
L(\lambda,\alpha,s)
=\sum_{0 \leq n \leq Y} \frac{e^{2\pi i \lambda n}}{(n+\alpha)^s}
+\delta_{\lambda} \frac{Y^{1-s}}{s-1}
+O_\lambda(Y^{-\sigma}),
\end{gather*}
where $\delta_{\lambda}=1$ if $\lambda=1$, $\delta_{\lambda}=0$ otherwise. 
\end{lemma}

\begin{proof}
This is \cite[Theorem 3.1.2]{LaurincikasGarunkstis2002} if $0<\lambda<1$. 
When $\lambda=1$, the Lerch zeta-function is the Hurwitz zeta-function, and hence the result is deduced from \cite[Theorem 3.1.3]{LaurincikasGarunkstis2002}. 
\end{proof}

\begin{proof}[Proof of Proposition \ref{p3.3}]
Note that we may assume that $0<\lambda \leq1$ without loss of generality. 
By the definition of $\psi_z(w)$, we have $|\psi_z(w)|=1$ and $|\psi_z(w)-\psi_z(w')| \leq |z||w-w'|$ for any $z, w, w' \in \mathbb{C}$. 
Hence we have 
\begin{gather*}
|E_1|
\leq \frac{2\pi}{T}
+\frac{|z|}{T} \int_{2\pi}^{T} 
\left|L(\lambda,\alpha,\sigma+it) - \sum_{0 \leq n \leq X} \frac{e^{2\pi i \lambda n}}{(n+\alpha)^{\sigma+it}}\right| \,dt.
\end{gather*}
Let $Y=T/\lambda$. 
Then we have $X<Y$ due to $\theta<2/3$. 
We deduce from Lemma \ref{l3.4} that 
\begin{gather*}
\left|L(\lambda,\alpha,\sigma+it)-\sum_{0 \leq n \leq X} \frac{e^{2\pi i \lambda n}}{(n+\alpha)^{\sigma+it}}\right|
\ll_\lambda \left|\sum_{X<n \leq Y} \frac{e^{2\pi i \lambda n}}{(n+\alpha)^{\sigma+it}}\right|
+\frac{T^{1-\sigma}}{t}
+T^{-\sigma}
\end{gather*}
for $\sigma>1/2$ and $2\pi \leq t \leq T$. 
Applying Cauchy's inequality, we obtain
\begin{equation}\label{e3.2}
E_1
\ll_\lambda \frac{1}{T}
+\frac{|z|}{T^{1/2}} 
\left( \int_{0}^{T} \left|\sum_{X<n \leq Y} \frac{e^{2\pi i \lambda n}}{(n+\alpha)^{\sigma+it}}\right|^2 \,dt \right)^{1/2}
+|z| T^{-\sigma} \log{T}
\end{equation}
By the inequality \cite[Theorem 2]{MontgomeryVaughan1974}, we see that 
\begin{gather*}
\int_{0}^{T} \left| \sum_{X<n \leq Y} \frac{e^{2\pi i \lambda n}}{(n+\alpha)^{\sigma+it}} \right|^2 \,dt
=\sum_{X<n \leq Y} \frac{T+O(n)}{(n+\alpha)^{2\sigma}}. 
\end{gather*}
Hence we have 
\begin{gather*}
\int_{0}^{T} \left| \sum_{X<n \leq Y} \frac{e^{2\pi i \lambda n}}{(n+\alpha)^{\sigma+it}} \right|^2 \,dt
\ll_\lambda T \sum_{n>X}(n+\alpha)^{-2\sigma}
\ll_{\sigma_1} \frac{T}{2\sigma-1} X^{1-2\sigma}. 
\end{gather*}
Therefore, by \eqref{e3.2}, we obtain 
\begin{gather*}
E_1
\ll_{\lambda,\sigma_1} \frac{1}{T}
+\frac{|z|}{\sqrt{2\sigma-1}} X^{1/2-\sigma}
+|z| T^{-\sigma} \log{T}.
\end{gather*}
Since $X=\exp((\log{T})^{3\theta/2})$, $|z| \leq2(\log{T})^\delta$, and $\sigma \geq 1/2+(\log{T})^{-\theta}$, we conclude 
\begin{gather*}
E_1 
\ll_{\lambda,\sigma_1} \exp\left(-\frac{1}{2} (\log{T})^{\theta/2}\right)
\end{gather*}
for any $T \geq T_0$ with a constant $T_0=T_0(\theta)>0$. 
\end{proof}

\begin{proposition}\label{p3.5}
Let $\lambda \in \mathbb{R}$ and $\alpha \in \mathfrak{S}$. 
Let $\sigma_1$ be a large fixed positive constant, and let $\theta, \delta>0$ with $\theta+\delta<1/4$. 
Then there exists $T_0=T_0(\alpha, \sigma_1, \theta, \delta)>0$ such that for all $U \geq T \geq T_0$ and for all $\sigma \in [1/2+(\log{T})^{-\theta}, \sigma_1]$, we have 
\begin{gather*}
\frac{1}{T} \int_{0}^{T} \psi_z\left(\sum_{0\leq n \leq X} \frac{e^{2\pi i \lambda n}}{(n+\alpha)^{\sigma+it}}\right) \,dt
=\frac{1}{U} \int_{0}^{U} \psi_z\left(\sum_{0\leq n \leq X} \frac{e^{2\pi i \lambda n}}{(n+\alpha)^{\sigma+iu}}\right) \,du
+E_2
\end{gather*}
for any $z \in \Omega$, where $X=\exp((\log{T})^{3\theta/2})$, and $E_2$ is estimated as 
\begin{gather*}
E_2
\ll\exp\left(-\frac{1}{2}(\log{T})^{\theta/2}\right).
\end{gather*}
\end{proposition}

\begin{proof}
For simplicity, we write
\begin{gather*}
F(t)
=F(t,\sigma,X;\lambda,\alpha)
=\sum_{0 \leq n \leq X} \frac{e^{2\pi i \lambda n}}{(n+\alpha)^{\sigma+it}}. 
\end{gather*}
Then we have 
\begin{align*}
\psi_z(F(t))
&=\exp\left( \frac{i \overline{z}}{2} F(t) + \frac{iz}{2} \overline{F(t)} \right)\\
&=\sum_{0 \leq \mu+\nu<N} \frac{1}{\mu! \nu!} 
\left(\frac{i \overline{z}}{2}\right)^\mu \left(\frac{iz}{2}\right)^\nu F(t)^\mu (\overline{F(t)})^\nu
+O\left(\frac{|z|^N}{N!} |F(t)|^N\right)
\end{align*}
for a sufficiently large positive even integer $N$. 
Hence $E_2$ is estimated as
\begin{align}\label{e3.3}
E_2
&=\sum_{0 \leq \mu+\nu<N} \frac{1}{\mu! \nu!} 
\left(\frac{i \overline{z}}{2}\right)^\mu \left(\frac{iz}{2}\right)^\nu H(\mu,\nu) \\
&\qquad
+O\left( \frac{|z|^N}{N!} \left(\frac{1}{T} \int_{0}^{T} |F(t)|^N \,dt
+\frac{1}{U} \int_{0}^{U} |F(u)|^N \,du\right) \right), \nonumber
\end{align}
where 
\begin{gather*}
H(\mu,\nu)
=\frac{1}{T} \int_{0}^{T} F(t)^\mu (\overline{F(t)})^\nu \,dt
-\frac{1}{U} \int_{0}^{U} F(u)^\mu (\overline{F(u)})^\nu \,du.
\end{gather*}
First, we estimate the first term of \eqref{e3.3}. 
We have 
\begin{align*}
H(\mu,\nu)
&\ll \frac{1}{T}
\mathop{ \sum_{0 \leq m_1 \leq X} \cdots \sum_{0 \leq m_\mu \leq X} 
\sum_{0 \leq n_1 \leq X} \cdots \sum_{0 \leq n_\nu \leq X} } 
\limits_{(m_1+\alpha) \cdots (m_\mu+\alpha) \neq (n_1+\alpha) \cdots (n_\nu+\alpha)}
\frac{1}{(m_1+\alpha)^\sigma \cdots (m_\mu+\alpha)^\sigma} \\
&\qquad 
\times \frac{1}{(n_1+\alpha)^\sigma \cdots (n_\nu+\alpha)^\sigma}
\left|\log\frac{(m_1+\alpha) \cdots (m_\mu+\alpha)}{(n_1+\alpha) \cdots (n_\nu+\alpha)}\right|^{-1}
\end{align*}
for $(\mu, \nu) \neq (0,0)$.
By condition \eqref{c2} of Definition \ref{d1.6}, this is 
\begin{gather*}
\ll \frac{1}{T} \left\{ \sum_{0 \leq n \leq X} \frac{1}{(n+\alpha)^\sigma} \right\}^N X^{\Omega(\alpha) N^2}
\leq \frac{1}{T} X^{(\Omega(\alpha)+1) N^2}
\end{gather*}
for $X \geq X(\alpha,\sigma_1)$ with a positive constant $X(\alpha,\sigma_1)$. 
By this and $H(0,0)=0$, the first term of \eqref{e3.3} is estimated as 
\begin{equation}\label{e3.4}
\ll \frac{(1+|z|^2)^{N/2}}{T} X^{(\Omega(\alpha)+1) N^2}.
\end{equation}
We consider the second term of \eqref{e3.3}. 
Let $N=2M$. 
Then we have 
\begin{align}\label{e3.5}
&\frac{1}{U} \int_{0}^{U} |F(u)|^{2M} \,du \\
&=\mathop{ \sum_{0 \leq m_1 \leq X} \cdots \sum_{0 \leq m_M \leq X}
\sum_{0 \leq n_1 \leq X} \cdots \sum_{0 \leq n_M \leq X} } 
\limits_{(m_1+\alpha) \cdots (m_M+\alpha) = (n_1+\alpha) \cdots (n_M+\alpha)}
\frac{e^{2\pi i \lambda m_1} \cdots e^{2\pi i \lambda m_M} e^{-2\pi i \lambda n_1} \cdots e^{-2\pi i \lambda n_M}}
{\{(n_1+\alpha) \cdots (n_M+\alpha)\}^{2\sigma}} \nonumber\\
&\qquad
+O\Bigg( \frac{1}{U}
\mathop{\sum_{0 \leq m_1 \leq X} \cdots \sum_{0 \leq m_M \leq X}
\sum_{0 \leq n_1\leq X} \cdots \sum_{0 \leq n_M \leq X} } 
\limits_{(m_1+\alpha) \cdots (m_M+\alpha) \neq (n_1+\alpha) \cdots (n_M+\alpha)}
\frac{1}{(m_1+\alpha)^\sigma \cdots (m_M+\alpha)^\sigma} \nonumber\\
&\qquad\qquad
\times\frac{1}{(n_1+\alpha)^\sigma \cdots (n_M+\alpha)^\sigma}
\left|\log \frac{(m_1+\alpha) \cdots (m_M+\alpha)}{(n_1+\alpha) \cdots (n_M+\alpha)}\right|^{-1} \Bigg). \nonumber
\end{align}
Due to condition \eqref{c1} of Definition \ref{d1.6}, the equation $(m_1+\alpha) \cdots (m_M+\alpha)=(n_1+\alpha) \cdots (n_M+\alpha)$ is equivalent to $\{ m_1,\ldots,m_M \}=\{ n_1,\ldots,n_M \}$.
Thus the diagonal term of \eqref{e3.5} is 
\begin{equation}\label{e3.6}
\ll M! \left(\sum_{0 \leq n \leq X} \frac{1}{(n+\alpha)^{2\sigma}}\right)^M
\leq M! \left(\frac{c(\alpha,\sigma_1)}{2\sigma-1}\right)^M
\end{equation}
with a positive constant $c(\alpha,\sigma_1)$.
The off-diagonal term is estimated as 
\begin{equation}\label{e3.7}
\ll \frac{1}{U} \left(\sum_{0 \leq n \leq X} \frac{1}{(n+\alpha)^{\sigma}}\right)^N X^{\Omega(\alpha) N^2}
\ll \frac{1}{T} X^{(\Omega(\alpha)+1) N^2}
\end{equation}
for $X \geq X(\alpha,\sigma_1)$ by applying condition \eqref{c2} again.
Therefore we have 
\begin{gather*}
E_2
\ll \frac{(1+|z|^2)^{N/2}}{T} X^{(\Omega(\alpha)+1) N^2}
+|z|^N \frac{(\frac{N}{2})!}{N!} \left(\frac{c(\alpha,\sigma_1)}{2\sigma-1}\right)^{N/2}
\end{gather*}
by \eqref{e3.4}, \eqref{e3.6}, and \eqref{e3.7}. 
By the assumption $\theta+\delta<1/4$, we have $\theta+2\delta<1/2-(3/4)\theta$.  
We take $\eta=\{(\theta+2\delta)+(1/2-(3/4)\theta)\}/2$ and $N=2\lfloor (\log{T})^{\eta} \rfloor$. 
Recalling $X=\exp((\log{T})^{3\theta/2})$, $|z| \leq2(\log{T})^\delta$, and $\sigma \geq 1/2+(\log{T})^{-\theta}$, we obtain 
\begin{gather*}
E_2 
\ll \exp\left(-\frac{1}{2} (\log{T})^{\theta/2}\right)
\end{gather*}
for $T \geq T_0$ with a positive constant $T_0=T_0(\alpha, \sigma_1, \theta, \delta)$. 
\end{proof}

\begin{proposition}\label{p3.6}
Let $\lambda \in \mathbb{R}$ and $\alpha \in \mathfrak{S}$. 
Then, for any $\sigma>1/2$ and $X \geq1$, we have 
\begin{gather*}
\lim_{U \to\infty} \frac{1}{U} \int_{0}^{U} 
\psi_z\left( \sum_{0 \leq n \leq X} \frac{e^{2\pi i \lambda n}}{(n+\alpha)^{\sigma+iu}} \right) \,du
=\prod_{0 \leq n \leq X} J_0(|z|(n+\alpha)^{-\sigma}).
\end{gather*}
\end{proposition}

\begin{proof}
Let $\gamma_n=-(2\pi)^{-1} \log(n+\alpha)$ for $n \geq0$.
Then $\{\gamma_n\}$ is linearly independent over $\mathbb{Q}$ by condition \eqref{c1} of Definition \ref{d1.6}. 
Then, by applying \cite[Lemma 2]{HeathBrown1992}, we have
\begin{align*}
\frac{1}{U} \int_{0}^{U} \psi_z\left(\sum_{0 \leq n \leq X} \frac{e^{2\pi i \lambda n}}{(n+\alpha)^{\sigma+iu}}\right) \,du
&=\frac{1}{U} \int_{0}^{U} \prod_{0 \leq n \leq X}
\psi_z\left(\frac{e^{2\pi i \lambda n}}{(n+\alpha)^{\sigma}} e^{2\pi i \gamma_n u}\right) \,du \\
&\to \prod_{0 \leq n \leq X} 
\int_{0}^{1} \psi_z\left(\frac{e^{2\pi i \lambda n}}{(n+\alpha)^{\sigma}} e^{2\pi i\theta}\right) \,d\theta 
\end{align*}
as $U \to\infty$.
Moreover, we see that 
\begin{gather*}
\int_{0}^{1} \psi_z\left(\frac{e^{2\pi i \lambda n}}{(n+\alpha)^{\sigma}} e^{2\pi i\theta}\right) \,d\theta
=\int_{0}^{1} \exp(i |z| (n+\alpha)^{-\sigma} \cos(2 \pi \varphi)) \,d\varphi 
\end{gather*}
by the changes of integral variables. 
The right-hand side is equal to $J_0(|z|(n+\alpha)^{-\sigma})$, and hence the result follows. 
\end{proof}

\begin{proposition}\label{p3.7}
Let $\lambda \in \mathbb{R}$ and $0<\alpha \leq1$. 
Let $\sigma_1$ be a large fixed positive constant, and let $\theta, \delta>0$. 
Then there exists $T_0=T_0(\theta)>0$ such that for all $T \geq T_0$ and for all $\sigma \in [1/2+(\log{T})^{-\theta},  \sigma_1]$, we have 
\begin{gather*}
\prod_{0 \leq n \leq X} J_0(|z| (n+\alpha)^{-\sigma})
=\widetilde{M}_\sigma(z; \alpha)
+E_3
\end{gather*}
for any $z \in \Omega$, where $X=\exp((\log{T})^{3\theta/2})$, and $E_3$ is estimated as 
\begin{gather*}
E_3
\ll_{\sigma_1} \exp\left(-\frac{1}{2} (\log{T})^{\theta/2}\right).
\end{gather*}
\end{proposition}

\begin{proof}[Proof of Theorem \ref{t3.1}]
It is easily directly deduced by combining Propositions \ref{p3.3}, \ref{p3.5}, \ref{p3.6}, and \ref{p3.7}. 
\end{proof}

\section{Proof of Theorem \ref{t1.7}}\label{s4}
We deduce Theorem \ref{t1.7} from Theorem \ref{t3.2} by approximating the characteristic function $1_R(z)$ by functions in the class $\Lambda$. 
For this, we use the following function $F_{a,b}(x)$. 

\begin{lemma}[Lemma 4.1 of \cite{Lester2014a}]\label{l4.1}
Let 
\begin{gather*}
K(x)
=\left(\frac{\sin \pi x}{\pi x}\right)^2.
\end{gather*}
Then, for any $a, b, \omega \in \mathbb{R}$ with $a<b$ and $\omega>0$, there exists a continuous function $F_{a,b,\omega}: \mathbb{R} \to \mathbb{R}$ such that the following conditions hold: 
\begin{enumerate}
\item[$\mathrm{(1)}$] 
$F_{a,b,\omega}(x)-1_{[a,b]}(x) \ll K(\omega(x-a))+K(\omega(x-b))$ for any $x \in \mathbb{R}$, 
\item[$\mathrm{(2)}$] 
$\displaystyle{ \int_{\mathbb{R}} (F_{a,b,\omega}(x)-1_{[a,b]}(x)) \,dx \ll \omega^{-1}}$, 
\item[$\mathrm{(3)}$] 
if $|x| \geq \omega$, then $\widetilde{F}_{a,b,\omega}(x)=0$,
\item[$\mathrm{(4)}$] 
$\widetilde{F}_{a,b,\omega}(x) \ll (b-a)+\omega^{-1}$. 
\end{enumerate}
Here 
\begin{gather*}
\widetilde{F}_{a,b,\omega}(x)
=\int_{\mathbb{R}} F_{a,b,\omega}(u) e^{ixu} \,|du|
\end{gather*}
is the Fourier transformation of $F_{a,b,\omega}(x)$ with $|du|=(2\pi)^{-1/2}du$. 
\end{lemma}

The function $F_{a,b,\omega}(x)$ is constructed by using the (refined) Beurding--Selberg functions 
\begin{gather*}
H(x)
=\left(\frac{\sin \pi x}{\pi}\right)^2 \left\{\sum_{m=-\infty}^{\infty} \sgn(m) (x-m)^{-2}+2 x^{-1}\right\}
\end{gather*}
which approximate well the signum function $\sgn(x)$; see \cite{Vaaler1985}. 

\begin{proof}[Proof of Theorem \ref{t1.7}]
We take the constants $\lambda$, $\alpha$, $\sigma_1$, $\theta$, $\delta>0$ as in the assumption of Theorem \ref{t3.2}. 
For $T \geq T_0(\alpha, \sigma_1, \theta, \delta)$, let $R$ be any rectangle 
\begin{gather*}
R
=\{ z=x+iy \mid a \leq x \leq b,~ c \leq y \leq d \}
\end{gather*}
with $b-a,~ d-c \gg (\log{T})^{-\delta}$.
We also take $\omega=(2\pi)^{-1}(\log{T})^\delta$. 
Note that condition $\mathrm{(4)}$ of Lemma \ref{l4.1} deduces $\widetilde{F}_{a,b,\omega}(x) \ll (b-a)$ in this case. 
Then, we define 
\begin{gather*}
\Phi(z)
=F_{a,b,\omega}(x) F_{c,d,\omega}(y).
\end{gather*}
Since we have 
\begin{gather*}
\Lambda
=\{ f \in L^1(\mathbb{C}) \mid \text{$f$ is continuous and $\widetilde{f} \in L^1(\mathbb{C})$} \}
\end{gather*}
by \cite[Section 9]{IharaMatsumoto2011a}, Lemma \ref{l4.1} deduces that $\Phi$ is a member of the class $\Lambda$. 
Moreover, we have 
\begin{equation}\label{e4.1}
\Phi(z)-1_R(z)
\ll K(\omega(x-a))+K(\omega(x-b))+K(\omega(y-c))+K(\omega(y-d))
\end{equation}
and $\widetilde{\Phi}(z) \ll (b-a)(d-c) = \mu_2(R)$ also by Lemma \ref{l4.1}. 

In the above setting, we apply Theorem \ref{t3.2}. 
Let $\epsilon_1>0$ arbitrarily, and we take $T_1=T_1(\alpha, \sigma_1, \theta, \delta, \epsilon_1) \geq T_0(\alpha, \sigma_1, \theta, \delta)$ so that $(\log{T_1})^{-\theta}<\epsilon_1$ holds. 
Then, for all $T \geq T_1$ and for all $1/2+\epsilon_1 \leq \sigma \leq \sigma_1$, we have 
\begin{equation}\label{e4.2}
\frac{1}{T} \int_{0}^{T} \Phi(L(\lambda, \alpha, \sigma+it)) \,dt
=\int_{\mathbb{C}} \Phi(z) M_\sigma(z; \alpha) \,|dz|
+E,
\end{equation}
where 
\begin{gather*}
E
\ll_{\lambda, \sigma_1} \exp\left(-\frac{1}{2} (\log{T})^{\theta/2}\right)
\int_{\Omega} \left|\widetilde{\Phi}(z)\right| \,|dz|
+\int_{\mathbb{C} \setminus \Omega} \left|\widetilde{\Phi}(z)\right| \,|dz|. 
\end{gather*}
By Lemma \ref{l4.1}, we see that  
\begin{equation}\label{e4.3}
E
\ll_{\lambda, \sigma_1} \exp\left(-\frac{1}{2}(\log{T})^{\theta/2}\right) (\log{T})^{2\delta} \mu_2(R) 
+0
\ll \mu_2(R)(\log{T})^{-\delta}
\end{equation}
for $T \geq T_2$ with a constant $T_2=T_2(\alpha, \sigma_1, \theta, \delta, \epsilon_1) \geq T_1$.
We estimate the error between the left-hand side of \eqref{e4.2} and 
\begin{gather*}
P_{\sigma, T}(R; \lambda, \alpha)
=\frac{1}{T} \int_{0}^{T} 1_R(L(\lambda, \alpha, \sigma+it)) \,dt.
\end{gather*}
Applying inequality \eqref{e4.1}, it is 
\begin{align*}
&\ll\frac{1}{T} \int_{0}^{T} K\left(\omega(\RE L(\lambda,\alpha,\sigma+it)-a)\right) \,dt \\
&\qquad
+\frac{1}{T} \int_{0}^{T} K\left(\omega(\RE L(\lambda,\alpha,\sigma+it)-b)\right) \,dt \\
&\qquad\quad
+\frac{1}{T} \int_{0}^{T} K\left(\omega(\IM L(\lambda,\alpha,\sigma+it)-c)\right) \,dt \\
&\qquad\qquad
+\frac{1}{T} \int_{0}^{T} K\left(\omega(\IM L(\lambda,\alpha,\sigma+it)-d)\right) \,dt. 
\end{align*}
Here, we consider only the term 
\begin{equation}\label{e4.4}
\frac{1}{T} \int_{0}^{T} K\left(\omega(\RE L(\lambda,\alpha,\sigma+it)-a)\right) \,dt,
\end{equation}
since the other terms are estimated in a similar way. 
We have 
\begin{gather*}
K(\omega x)
=\frac{2}{\omega^2} \int_{0}^{\omega} (\omega-u) \cos(2\pi xu) \,du
=\frac{2}{\omega^2} \RE \int_{0}^{\omega} (\omega-u) e^{2\pi i xu} \,du,
\end{gather*}
and hence \eqref{e4.4} is estimated as 
\begin{align*}
&\ll \frac{1}{\omega^2} \int_{0}^{\omega}(\omega-u)
\left|\frac{1}{T} \int_{0}^{T} \exp(2\pi iu(\RE L(\lambda,\alpha,\sigma+it)-a)) \,dt\right| \,du\\
&=\frac{1}{\omega^2} \int_{0}^{\omega} (\omega-u)
\left|\frac{1}{T} \int_{0}^{T} \exp(2\pi iu \RE L(\lambda,\alpha,\sigma+it)) \,dt\right| \,du\\
&\ll_{\lambda, \sigma_1} \frac{1}{\omega^2} \int_{0}^{\omega} (\omega-u) 
\left\{ \left|\widetilde{M}_\sigma(2\pi u;\alpha)\right|+\exp\left(-\frac{1}{2}(\log{T})^{\theta/2}\right) \right\} \,du\\
&\ll_{\alpha, \epsilon_1} \frac{1}{\omega} 
+\exp\left(-\frac{1}{2}(\log{T})^{\theta/2}\right)
\ll(\log{T})^{-\delta}. 
\end{align*}
Here we used Theorem \ref{t3.1} for the second inequality and used the fact that $\widetilde{M}_\sigma(z;\alpha)$ rapidly decreases for the third inequality. 
Therefore we have 
\begin{equation}\label{e4.5}
\frac{1}{T} \int_{0}^{T} \Phi(L(\lambda, \alpha, \sigma+it)) \,dt - P_{\sigma, T}(R; \lambda, \alpha)
\ll_{\lambda,\alpha,\sigma_1,\epsilon_1} (\log{T})^{-\delta}. 
\end{equation}
At the end of the proof, we estimate the difference
\begin{gather*}
\int_{\mathbb{C}} \Phi(z) M_\sigma(z; \alpha) \,|dz|
-\int_{R} M_\sigma(z; \alpha) \,|dz|
=\int_{\mathbb{C}} (\Phi(z)-1_R(z)) M_\sigma(z; \alpha) \,|dz|
\end{gather*}
coming from the right-hand side by applying \eqref{e4.1} again. 
We have 
\begin{equation}\label{e4.6}
\int_{\mathbb{C}} K(\omega(x-a)) M_\sigma(z;\alpha) \,|dz|
=\int_{\mathbb{R}} K(\omega(x-a)) m_\sigma(x;\alpha) \,|dx|, 
\end{equation}
where 
\begin{gather*}
m_\sigma(x;\alpha)
=\int_{\mathbb{R}} M_\sigma(x+iy;\alpha) \,|dy|.
\end{gather*}
By estimate \eqref{e2.5}, the function $m_\sigma(x;\alpha)$ is estimated as 
\begin{gather*}
m_\sigma(x;\alpha)
\ll_{\alpha,\epsilon_1} 1. 
\end{gather*}
Hence we obtain that \eqref{e4.6} is 
\begin{gather*}
\ll_{\alpha,\epsilon_1} \int_{\mathbb{R}} K(\omega(x-a)) \,|dx|
\ll \omega^{-1}
=(\log{T})^{-\delta}.
\end{gather*}
The remaining terms are estimated similarly. 
We have 
\begin{equation}\label{e4.7}
\int_{\mathbb{C}} \Phi(z) M_\sigma(z; \alpha) \,|dz|
-\int_{R} M_\sigma(z; \alpha) \,|dz|
\ll_{\alpha,\epsilon_1} (\log{T})^{-\delta}.
\end{equation}

Therefore, we conclude 
\begin{gather*}
P_{\sigma, T}(R; \lambda, \alpha)
=\int_{R} M_\sigma(z; \alpha) \,|dz|
+O_{\lambda,\alpha,\sigma_1,\epsilon_1} \left((\mu_2(R)+1) (\log{T})^{-\delta}\right)
\end{gather*}
by estimates \eqref{e4.3}, \eqref{e4.5}, and \eqref{e4.7}. 
We take $\theta=\epsilon/2$ and $\delta=1/4-\epsilon$ for each $0<\epsilon<1/4$. 
Then we have $\theta,\delta>0$ and $\theta+\delta<1/4$, and the desired result follows. 
\end{proof}

\section{Further results on the density function}\label{s5}
We have obtained several results on the function $M_\sigma(z; \alpha)$ in Section \ref{s2}. 
In this section, we study it as a function of $\sigma$. 
The ultimate goal of this section is the following estimate. 

\begin{proposition}\label{p5.1}
Let $0<\alpha \leq1$ and $\sigma>1/2$. 
Then there exists a positive constant $c(\sigma)$ such that for any integers $k, l, m \geq0$, we have 
\begin{gather*}
\frac{\partial^{k+l+m}}{\partial x^k \partial y^l \partial \sigma^m} \widetilde{M}_\sigma(z; \alpha)
\ll_{k,l,m} \exp\left(-c(\sigma) |z|^{1/\sigma}\right),
\qquad 
|z| \to\infty
\end{gather*}
with $z=x+iy$. 
\end{proposition}

To begin with, we recall that   
\begin{align*}
J_0(|z|(n+\alpha)^\sigma)
&=\int_{0}^{1} \psi_z((n+\alpha)^{-\sigma} e^{2\pi i\theta}) \,d\theta\\
&=\int_{0}^{1} \exp\{ix(n+\alpha)^{-\sigma} \cos(2\pi \theta)+iy(n+\alpha)^{-\sigma} \sin(2\pi \theta)\} \,d\theta.
\end{align*}
We then define 
\begin{align}\label{e5.1}
&\widetilde{M}_n(s,z_1,z_2; \alpha) \\
&=\int_{0}^{1} \exp\{i z_1(n+\alpha)^{-s} \cos(2\pi \theta)+i z_2(n+\alpha)^{-s} \sin(2\pi \theta)\} \,d\theta \nonumber\\
&=1
-\frac{z_1^2+z_2^2}{4} (n+\alpha)^{-2s} \nonumber\\
&\qquad\quad
+\int_{0}^{1} \sum_{k=3}^{\infty} \frac{i^k}{k!}(n+\alpha)^{-ks} 
(z_1 \cos(2\pi \theta)+z_2 \sin(2\pi \theta))^k \,d\theta \nonumber
\end{align}
for $s, z_1, z_2 \in \mathbb{C}$ with $\RE s>0$. 
Note that $\widetilde{M}_n(\sigma,x,y; \alpha)=J_0(|x+iy|(n+\alpha)^\sigma)$ if $\sigma, x, y \in \mathbb{R}$. 
We investigate the function 
\begin{equation}\label{e5.2}
\widetilde{M}(s,z_1,z_2; \alpha)
=\prod_{n=0}^{\infty} \widetilde{M}_n(s,z_1,z_2; \alpha).
\end{equation}

\begin{lemma}\label{l5.2}
If we fix two of the variables, the local parts $\widetilde{M}_n(s,z_1,z_2;\alpha)$ are holomorphic with respect to the remaining variable for any $s, z_1, z_2 \in \mathbb{C}$ with $\RE s>0$. 
Let $K$ be any compact subset on the half plane $\{\RE s>1/2\}$, and let $K_1, K_2$ be any compact subsets on $\mathbb{C}$. 
Then the infinite product \eqref{e5.2} uniformly converges on $K \times K_1 \times K_2$. 
Therefore, if we again fix two of the variables, the function $\widetilde{M}(s,z_1,z_2;\alpha)$ is holomorphic with respect to the remaining variable for any $s, z_1, z_2 \in \mathbb{C}$ with $\RE s>1/2$. 
\end{lemma}

\begin{proof}
The first statement is clear from the definition \eqref{e5.1}. 
Assume that $(s, w_1, w_2)$ varies in $K \times K_1 \times K_2$, and let $\sigma_0$ be the smallest real parts of $s \in K$. 
Then there exists a sufficiently large constant $N=N(K,K_1,K_2;\alpha)$ such that 
\begin{gather*}
\widetilde{M}_n(s,z_1,z_2; \alpha)
=1+O_{K_1,K_2}\left((n+\alpha)^{-2\sigma_0}\right)
\end{gather*}
for $n \geq N$ by \eqref{e5.1}. 
Since the series $\sum_{n} (n+\alpha)^{-2\sigma_0}$ converges due to $\sigma_0>1/2$, the second statement follows. 
The last statement directly follows from the preceding two statement. 
\end{proof}

In the remaining part of this section, we prove the following estimate. 

\begin{lemma}\label{l5.3}
Let $0<\alpha \leq1$ and $\sigma>1/2$. 
Then there exist positive constants $K(\sigma;\alpha)$ and $c(\sigma)$ such that for any $x,y \in \mathbb{R}$ with $|x|+|y|\geq K(\sigma;\alpha)$, we have
\begin{gather*}
\left|\widetilde{M}(s,z_1,z_2;\alpha)\right|
\leq \exp\left(-c(\sigma)(|x|+|y|)^{1/\RE(s)}\right)
\end{gather*}
for any $s, z_1, z_2 \in \mathbb{C}$ with $|s-\sigma|<1/(x^2+y^2)$, $|z_1-x|<1/2$, $|z_2-y|<1/2$. 
\end{lemma}

\begin{proof}
Let $n_0$ be a real number satisfying 
\begin{gather*}
n_0+\alpha
=\left(\frac{|x|+|y|}{c_0}\right)^{1/\RE(s)},
\end{gather*}
where $c_0$ is a positive absolute constant suitably chosen later. 
At first, we take $c_0<1/2$ and $K(\sigma;\alpha) \geq1$. 
Then, for $n \geq n_0$, we have 
\begin{gather*}
(|x|+|y|) (n+\alpha)^{-\RE(s)}
\leq(|x|+|y|)n_0^{-\RE(s)}
=c_0.
\end{gather*}
Hence we see that 
\begin{gather*}
\left|\frac{z_1^2+z_2^2}{4} (n+\alpha)^{-2s}\right|
\leq c_0^2, \\
\left|\int_{0}^{1} \sum_{k=3}^{\infty} 
\frac{i^k}{k!} (n+\alpha)^{-ks} (z_1 \cos(2\pi \theta)+z_2 \sin(2\pi \theta))^k \,d\theta\right|
\leq 4c_0^3 
\end{gather*}
for $n \geq n_0$ and $|z_1-x|<1/2$, $|z_2-y|<1/2$. 
By formula \eqref{e5.1}, we can define $\Log\widetilde{M}_n(s,z_1,z_2; \alpha)$ and obtain 
\begin{equation}\label{e5.3}
\Log\widetilde{M}_n(s,z_1,z_2; \alpha)
=-\frac{z_1^2+z_2^2}{4} (n+\alpha)^{-2s}
+O\left((|x|+|y|)^3 (n+\alpha)^{-3\RE(s)}\right).
\end{equation}
Since we see that  
\begin{align*}
&\left|\frac{z_1^2+z_2^2}{4} (n+\alpha)^{-2s}-\frac{x^2+y^2}{4} (n+\alpha)^{-2\RE(s)}\right|\\
&\leq\left|\frac{z_1^2+z_2^2}{4}-\frac{x^2+y^2}{4}\right| (n+\alpha)^{-2\RE(s)}
+\frac{x^2+y^2}{4} \left|(n+\alpha)^{-2s}-(n+\alpha)^{-2\RE(s)}\right|\\
&\ll(|x|+|y|) (n+\alpha)^{-2\RE(s)}
+\log(n+\alpha) (n+\alpha)^{-2\RE(s)}
\end{align*}
for $|s-\sigma|<1/(x^2+y^2)$, there exists an absolute constant $A>0$ such that 
\begin{align*}
&\left|\Log\widetilde{M}_n(s,z_1,z_2; \alpha)+\frac{x^2+y^2}{4} (n+\alpha)^{-2\RE(s)}\right|\\
&\leq A \left((|x|+|y|)(n+\alpha)^{-\RE(s)}+(|x|+|y|)^{-1}\right) (|x|+|y|)^2 (n+\alpha)^{-2\RE(s)}\\
&\qquad
+A \log(n+\alpha) (n+\alpha)^{-2\RE(s)}. 
\end{align*}
Then, assuming further $c_0<(32A)^{-1}$ and $K(\sigma,\alpha) \geq 32A$, we have  
\begin{gather*}
A\left((|x|+|y|) (n+\alpha)^{-\RE(s)} + (|x|+|y|)^{-1}\right)
\leq A (c_0+K^{-1})
\leq \frac{1}{16} 
\end{gather*}
for $|x|+|y| \geq K(\sigma;\alpha)$. 
Hence, by \eqref{e5.3}, we obtain 
\begin{gather*}
\RE \Log \widetilde{M}_n(s,z_1,z_2; \alpha)
\leq-\frac{(|x|+|y|)^2}{16} (n+\alpha)^{-2\RE(s)} 
+A \log(n+\alpha) (n+\alpha)^{-2\RE(s)}
\end{gather*}
for $n \geq n_0$, and moreover, 
\begin{align}\label{e5.4}
\left|\prod_{n \geq n_0} \widetilde{M}_n(s,z_1,z_2; \alpha)\right| 
&\leq \exp\Bigg(-\frac{(|x|+|y|)^2}{16} \sum_{n \geq n_0}(n+\alpha)^{-2\RE(s)} \\
&\qquad\qquad\qquad
+A \sum_{n \geq n_0} \log(n+\alpha) (n+\alpha)^{-2\RE(s)}\Bigg). \nonumber
\end{align}
Since we have 
\begin{gather*}
|\RE(s)-\sigma|
\leq |s-\sigma|
<\frac{1}{x^2+y^2}
\leq \frac{4}{K(\sigma;\alpha)^2}, 
\end{gather*}
we deduce 
\begin{gather*}
2\RE(s)-1
>\sigma-1/2
>0
\end{gather*}
if we take $K(\sigma;\alpha)>4/\sqrt{2\sigma-1}$. 
Then, for the first term of \eqref{e5.4}, we have 
\begin{align*}
(|x|+|y|)^2 \sum_{n \geq n_0} (n+\alpha)^{-2\RE(s)}
&\geq \frac{1}{2\RE(s)-1} (n_0+\alpha)^{1-2\RE(s)} \\
&=\frac{c_0^{\frac{2\RE(s)-1}{\RE(s)}}}{2\RE(s)-1} (|x|+|y|)^{1/\RE(s)}. 
\end{align*}
On the other hand, we have 
\begin{align*}
&\sum_{n \geq n_0} \log(n+\alpha) (n+\alpha)^{-2\RE(s)} \\
&\leq \frac{1}{2\RE(s)-1} \log(n_0+\alpha) (n_0+\alpha)^{1-2\RE(s)} \\
&\qquad\qquad
\times \left(1+\frac{2\RE(s)-1}{n_0+\alpha}+\frac{1}{(2\RE(s)-1) \log(n_0+\alpha)}\right). 
\end{align*}
Therefore, we see that there exists $K(\sigma;\alpha)>0$ such that for any $|x|+|y| \geq K(\sigma;\alpha)$, 
\begin{gather*}
\sum_{n \geq n_0} \log(n+\alpha) (n+\alpha)^{-2\RE(s)}
\leq \Theta \frac{c_0^{\frac{2\RE(s)-1}{\RE(s)}}}{2\RE(s)-1} (|x|+|y|)^{1/\RE(s)}
\end{gather*}
with a positive constant $\Theta < (16A)^{-1}$. 
Thus \eqref{e5.4} deduces 
\begin{equation}\label{e5.5}
\left|\prod_{n \geq n_0} \widetilde{M}_n(s,z_1,z_2; \alpha)\right|
\leq \exp\left(-c(\sigma) (|x|+|y|)^{1/\RE(s)}\right),
\end{equation}
where $c(\sigma)$ is a positive constant depending only on $\sigma$. 

Then, we estimate $\widetilde{M}_n(s,z_1,z_2; \alpha)$ for $n<n_0$. 
By \eqref{e5.1}, we have 
\begin{align*}
\left|\widetilde{M}_n(s,z_1,z_2; \alpha)\right|
&\leq \exp\left\{ (|\IM(z_1)|+|\IM(z_2)|) (n+\alpha)^{-\RE(s)} \right\} \\
&\leq \exp\left((n+\alpha)^{-\RE(s)}\right)
\end{align*}
since $|z_1-x|<1/2$ and $|z_2-y|<1/2$ with $x, y \in \mathbb{R}$. 
Therefore the contribution from the terms for $n<n_0$ is estimated as 
\begin{align}\label{e5.6}
\left|\prod_{n<n_0} \widetilde{M}_n(s,z_1,z_2; \alpha)\right|
&\leq \exp\left(\sum_{n<n_0} (n+\alpha)^{-\RE(s)}\right) \nonumber\\
&\leq \exp\left(c'(\sigma;\alpha) (|x|+|y|)^{1/2\RE(s)}\right),
\end{align}
where $c'(\sigma;\alpha)$ is a suitable positive constant. 
Hence we have the desired result from estimates \eqref{e5.5} and \eqref{e5.6}.  
\end{proof}

\begin{proof}[Proof of Proposition \ref{p5.1}]
By Lemma \ref{l5.2}, we can apply Cauchy's integral formula for the function $\widetilde{M}(s,z_1,z_2;\alpha)$. 
Proposition \ref{p5.1} is easily deduced from this and the estimate of Lemma \ref{l5.3}. 
\end{proof}

\begin{corollary}\label{c5.4}
Let $0<\alpha \leq1$ and $\sigma>1/2$. 
Then, for any integer $m \geq0$, the function
\begin{equation}\label{e5.7}
\frac{\partial^{m}}{\partial \sigma^m} M_\sigma(z; \alpha)
\end{equation}
belongs to $\mathcal{S}$ as a function in $x$ and $y$ with $z=x+iy$. 
\end{corollary}

\begin{proof}
By Proposition \ref{p5.1}, the function
\begin{equation}\label{e5.8}
\frac{\partial^{m}}{\partial \sigma^m} \widetilde{M}_\sigma(z; \alpha)
\end{equation}
belongs to $\mathcal{S}$. 
Thus we have 
\begin{align*}
\frac{\partial^{m}}{\partial \sigma^m} M_\sigma(z; \alpha)
&=\frac{\partial^{m}}{\partial \sigma^m} \int_{\mathbb{C}} \widetilde{M}_\sigma(w; \alpha) \psi_{-z}(w) \,|dw| \\
&=\int_{\mathbb{C}} \frac{\partial^{m}}{\partial \sigma^m} \widetilde{M}_\sigma(w; \alpha) \psi_{-z}(w) \,|dw|.
\end{align*}
In other words, the Fourier inverse of function \eqref{e5.8} is equal to function \eqref{e5.7}. 
Hence function \eqref{e5.7} belongs to $\mathcal{S}$. 
\end{proof}

\section{Proof of Theorem \ref{t1.8}}\label{s6}
Finally we prove Theorem \ref{t1.8}. 
The following proposition is an analogue of \cite[Theorem 1.1.3]{Guo1996b} for Lerch zeta-functions and is a key for the proof of Theorem \ref{t1.8}. 

\begin{proposition}\label{p6.1}
Let $\lambda \in \mathbb{R}$ and $\alpha \in \mathfrak{S}$. 
Let $\sigma_1$ be a large fixed positive constant. 
Let $\epsilon_1>0$ be a small fixed real number. 
Then there exists $T_0=T_0(\lambda, \alpha, \sigma_1, \epsilon_1)>0$ such that for all $T \geq T_0$ and for all $1/2+\epsilon_1 \leq \sigma \leq \sigma_1$, we have 
\begin{equation}\label{e6.1}
\frac{1}{T} \int_{0}^{T} \log|L(\lambda, \alpha, \sigma+it)| \,dt
=\int_{\mathbb{C}} \log|z| M_\sigma(z; \alpha) \,|dz|
+O_{\lambda,\alpha,\sigma_1,\epsilon_1} \left((\log{T})^{-A}\right)
\end{equation}
with an absolute constant $A>0$. 
\end{proposition}

Before the proof of Proposition \ref{p6.1}, we estimate the upper bound of the integral
\begin{gather*}
\int_{0}^{T} \log^2 |L(\lambda,\alpha,\sigma+it)| \,dt 
\end{gather*}
for $\sigma>1/2$ in Section \ref{s6.1}. 
It is used to bound the error term coming from the left-hand side of \eqref{e6.1}. 
We then prove Proposition \ref{p6.1} in Section \ref{s6.2}. 
Proposition \ref{p6.1} is connected to Theorem \ref{t1.8} by Lemma \ref{l6.8} which is a consequence of well-known Littlewood's lemma \cite[Lemma 8.4.9]{LaurincikasGarunkstis2002}; see also \cite[p.\,221]{Titchmarsh1986}. 
The proof of Theorem \ref{t1.8} is completed in Section \ref{s6.3}. 

\subsection{Mean square of the logarithm of the Lerch zeta-function}\label{s6.1}

\begin{proposition}\label{p6.2}
Let $\lambda \in \mathbb{R}$ and $0<\alpha \leq1$. 
Let $\sigma_1>0$ be a large fixed real number. 
Let $\epsilon_1>0$ be a small fixed real number. 
Then there exists $T_0=T_0(\lambda,\alpha,\sigma_1,\epsilon_1)>0$ such that for all $T \geq T_0$ and for all $1/2+\epsilon_1 \leq \sigma \leq \sigma_1$, we have  
\begin{gather*}
\int_{0}^{T} \log^2|L(\lambda,\alpha,\sigma+it)|\,dt 
\ll_{\lambda,\alpha,\sigma_1,\epsilon_1} T.
\end{gather*}
\end{proposition}

The proof is based on the method in \cite[Section 5]{LamzouriLesterRadziwill2019}. 
Let 
\begin{gather*}
f(z)
=\alpha^z L(\lambda,\alpha,z).
\end{gather*}
By \cite[Lemma 8.5.2]{LaurincikasGarunkstis2002}, the function $f(z)$ satisfies
\begin{equation}\label{e6.2}
1-\frac{\alpha}{x-1}
<|f(z)|
<1+\frac{\alpha}{x-1}, 
\qquad z=x+iy
\end{equation}
for $x>1+\alpha$. 
Moreover, let $\sigma_0=\max(\sigma_1,3)$, and let 
\begin{gather*}
r
=\sigma_0-\frac{1}{2} \left(\sigma+\frac{1}{2}\right), 
\quad
\delta
=\frac{1}{\sigma_1+4} \left(\sigma-\frac{1}{2}\right), 
\quad
r_1
=r-\delta,
\quad
r_2
=r-2\delta.
\end{gather*}
Note that we have $r>r_1>r_2>0$, $0<\delta<1$, and $r+2\delta<\sigma_0$. 
We define $s_0=s_0(t)=\sigma_0+it$. 
If $t \geq \sigma_0$, the function $f(z)$ is analytic in $|z-s_0| \leq r+2\delta$ and $1/2 \leq |f(s_0)| \leq 3/2$ by inequality \eqref{e6.2}. 
Then, for $t \geq \sigma_0$ and $0<u \leq r+2\delta$, we define  
\begin{gather*}
M_u(t)
=\max_{\substack{z \\ |z-s_0(t)| \leq u}} \left|\frac{f(z)}{f(s_0(t))}\right| 
+3
\end{gather*}
and
\begin{gather*}
N_u(t)
=\sum_{\substack{\rho \\ |\rho-s_0(t)| \leq u}} 1,
\end{gather*}
where $\rho$ runs through the zeros of $L(\lambda,\alpha,s)$.
We begin with the following formula. 

\begin{lemma}\label{l6.3}
Let $\lambda \in \mathbb{R}$ and $0<\alpha \leq1$. 
Let $\sigma_1>0$ be a large fixed real number. 
Let $\epsilon_1>0$ be a small fixed real number. 
Then for all $t \geq \sigma_0$ and for all $1/2+\epsilon_1 \leq \sigma \leq \sigma_1$, we have  
\begin{gather*}
\log|L(\lambda,\alpha,\sigma+it)|
=\sum_{|\rho-s_0(t)| \leq r_1} \log|\sigma+it-\rho|
+O_{\alpha,\sigma_1,\epsilon}\left( \log M_r(t) \right),
\end{gather*}
where $\rho$ runs through the zeros of $L(\lambda,\alpha,s)$.
\end{lemma}

\begin{proof}
We apply \cite[Lemma 2.2.1]{Guo1996b}. 
Then we have for $|z-s_0| \leq r_2$, 
\begin{align}\label{e6.3}
&\left|\frac{f'}{f}(z)-\sum_{|\rho-s_0| \leq r_1} \frac{1}{z-\rho}\right| \\
&\leq \frac{36r_1}{(r_1-r_2)^2} \left\{\log M_r(t)+N_{r_1}(t) \log \left(\frac{r_1}{r-r_1}+1\right)\right\} \nonumber\\
&\leq \frac{36\sigma_0}{\delta^2} \left(\log M_r(t)+N_{r_1}(t) \frac{\sigma_0}{\delta}\right). \nonumber
\end{align}
Jensen's formula yields that the equation
\begin{gather*}
\int_{0}^{r} \frac{N_x(t)}{x} \,dx
+\log|f(s_0)|
=\frac{1}{2\pi} \int_{0}^{2\pi} \log|f(s_0+re^{i \theta})| \,d\theta
\end{gather*}
holds, and its left-hand side is estimated as 
\begin{gather*}
\geq \int_{r_1}^{r} \frac{N_x(t)}{x} \,dx
+\log|f(s_0)|
\geq N_{r_1}(t) \frac{r-r_1}{r}
+\log|f(s_0)|.
\end{gather*}
Also, the right-hand side is 
\begin{gather*}
\leq \frac{1}{2\pi} \int_{0}^{2\pi} \log\left|\frac{f(s_0+re^{i \theta})}{f(s_0)}\right| \,d\theta
+\log|f(s_0)|
\leq \log M_r(t)
+\log|f(s_0)|.
\end{gather*}
Hence we have 
\begin{equation}\label{e6.4}
N_{r_1}(t)
\leq \frac{r}{r-r_1} \log M_r(t)
\leq \frac{\sigma_0}{\delta} \log M_r(t). 
\end{equation}
By \eqref{e6.3} and \eqref{e6.4}, we obtain 
\begin{gather*}
\frac{f'}{f}(z)-\sum_{|\rho-s_0| \leq r_1} \frac{1}{z-\rho}
\ll_{\sigma_1,\epsilon_1} \log M_r(t)
\end{gather*}
for all $|z-s_0| \leq r_2$, since $\sigma_0 \ll_{\sigma_1} 1$ and $\delta \gg_{\sigma_1,\epsilon_1} 1$. 
Note that $|\sigma+it-s_0|=\sigma_0-\sigma \leq r_2$. 
Therefore, integrating from $s_0$ to $\sigma+it$, we have 
\begin{align*}
&\log|f(\sigma+it)|-\sum_{|\rho-s_0| \leq r_1} \log|\sigma+it-\rho|\\
&=\log|f(s_0)|-\sum_{|\rho-s_0| \leq r_1} \log|s_0-\rho|
+\int_{s_0}^{\sigma+it} \left(\frac{f'}{f}(z)-\sum_{|\rho-s_0|\leq r_1} \frac{1}{z-\rho}\right) \,dz\\
&=\log|f(s_0)|-\sum_{|\rho-s_0| \leq r_1} \log|s_0-\rho|
+O_{\sigma_1,\epsilon_1}\left(\log M_r(t)\right). 
\end{align*}
Since $1/2 \leq |f(s_0)| \leq 3/2$, we have $\log|f(s_0)| \ll1$. 
Also, since $1 \leq |s_0-\rho| \leq \sigma_0$ for all zeros $\rho$ with $|s_0-\rho|<r_1$, we have 
\begin{gather*}
\sum_{|\rho-s_0| \leq r_1} \log|s_0-\rho|
\ll_{\sigma_1} N_{r_1}(t)
\ll_{\sigma_1,\epsilon_1} \log M_r(t)
\end{gather*}
by \eqref{e6.4}. 
By the definition of $f(z)$, the desired result follows. 
\end{proof}

By Lemma \ref{l6.3}, we have 
\begin{align}\label{e6.5}
&\int_{T}^{2T} \log^2|L(\lambda,\alpha,\sigma+it)| \,dt\\
&\ll \int_{T}^{2T} \left(\sum_{|\rho-s_0(t)| \leq r_1} \log|\sigma+it-\rho|\right)^2 \,dt
+\int_{T}^{2T} \log^2 M_r(t) \,dt \nonumber
\end{align}
for $T \geq \sigma_0$. 
Next we estimate two integrals of the right-hand side of \eqref{e6.5}.

\begin{lemma}\label{l6.4}
Let $\lambda \in \mathbb{R}$ and $0<\alpha \leq1$. 
Let $\sigma_1$ be a large fixed positive constant. 
Let $\epsilon_1>0$ be a small fixed real number. 
Then there exists $T_0=T_0(\lambda, \alpha, \sigma_1, \epsilon_1)>0$ such that for all $T \geq T_0$ and for all $1/2+\epsilon_1 \leq \sigma \leq \sigma_1$, we have 
\begin{equation}\label{e6.6}
\int_{T}^{2T} \left(\sum_{|\rho-s_0(t)| \leq r_1} \log|\sigma+it-\rho|\right)^2 \,dt 
\ll_{\lambda,\alpha,\sigma_1,\epsilon_1} T,
\end{equation}
where $\rho$ runs through the zeros of $L(\lambda,\alpha,s)$.
\end{lemma}

\begin{proof}
We have 
\begin{align*}
&\int_{T}^{2T} \left(\sum_{|\rho-s_0(t)| \leq r_1} \log|\sigma+it-\rho|\right)^2 \,dt \\
&\leq \sum_{n=\lfloor T\rfloor}^{\lfloor2T\rfloor+1} \int_{n}^{n+1} 
\left(\sum_{|\rho-s_0(t)| \leq r_1} \log|\sigma+it-\rho|\right)^2 \,dt. 
\end{align*}
According to the method in \cite{LamzouriLesterRadziwill2019}, let 
\begin{gather*}
\mathcal{D}_n
=\bigcup_{l=0}^{\lfloor1/\sqrt{\delta}\rfloor+1} D_r\left(\sigma_0+i(n+l\sqrt{\delta})\right),
\end{gather*}
where $D_r(c)=\{z \in \mathbb{C} \mid |z-c| \leq r\}$. 
Then we see that 
\begin{gather*}
\mathcal{D}_n \supset \bigcup_{n\leq t\leq n+1} D_{r_1}(s_0(t)),
\end{gather*}
and therefore, 
\begin{align*}
&\int_{n}^{n+1} \left(\sum_{|\rho-s_0(t)| \leq r_1} \log|\sigma+it-\rho|\right)^2 \,dt \\
&\leq \int_{n}^{n+1} \left(\sum_{\rho \in \mathcal{D}_n} \log|\sigma+it-\rho|\right)^2 \,dt\\
&\leq \left(\sum_{\rho \in \mathcal{D}_n} \left(\int_{n}^{n+1} \log^2|\sigma+it-\rho| \,dt\right)^{1/2} \right)^2
\end{align*}
due to Minkowski's inequality. 
If $n \leq t \leq n+1$ and $\rho=\beta+i \gamma \in \mathcal{D}_n$, we have 
\begin{gather*}
|t-\gamma|
\leq |\sigma+it-\rho|
\leq c
\end{gather*}
for a constant $c=c(\sigma_1)>1$. 
Thus 
\begin{gather*}
\log^2|\sigma+it-\rho|
\leq \log^2|t-\gamma|+\log^2{c}.
\end{gather*}
Moreover, we see that 
\begin{gather*}
\int_{n}^{n+1} \log^2|t-\gamma| \,dt
\leq \int_{n-r-\gamma}^{n+r+2-\gamma} \log^2{x} \,dx
\leq \int_{-2r-2}^{2r+2} \log^2{x} \,dx
\ll_{\sigma_1}1,
\end{gather*}
since $n-1 \leq \gamma \leq n+2+r$ for $\rho=\beta+i \gamma \in \mathcal{D}_n$. 
From the above, we deduce that the left-hand side of \eqref{e6.6} is 
\begin{align}\label{e6.7}
\int_{T}^{2T}\left(\sum_{|\rho-s_0(t)| \leq r_1} \log|\sigma+it-\rho|\right)^2 \,dt
&\ll_{\sigma_1} \sum_{n=\lfloor T\rfloor}^{\lfloor2T\rfloor+1} \left(\sum_{\rho \in \mathcal{D}_n} 1\right)^2 \\
&\leq \sum_{n=\lfloor T\rfloor}^{\lfloor2T\rfloor+1} 
\left(\sum_{l=0}^{\lfloor1/\sqrt{\delta}\rfloor+1} N_r\left(n+l\sqrt{\delta}\right)\right)^2. \nonumber
\end{align}
By an argument similar to \eqref{e6.4}, we see that 
\begin{gather*}
N_r\left(n+l\sqrt{\delta}\right)
\leq \frac{\sigma_0}{\delta} \log M_{r+\delta} \left(n+l\sqrt{\delta}\right).
\end{gather*}
Applying this inequality, we obtain 
\begin{equation}\label{e6.8}
\sum_{n=\lfloor T\rfloor}^{\lfloor2T\rfloor+1} 
\left(\sum_{l=0}^{\lfloor1/\sqrt{\delta}\rfloor+1} N_r\left(n+l\sqrt{\delta}\right)\right)^2
\ll_{\sigma_1,\epsilon_1} \sum_{n=\lfloor T\rfloor}^{\lfloor2T\rfloor+1} \log^2\left(2 \max_{z \in E_n}|f(z)|+3\right),
\end{equation}
where 
\begin{gather*}
E_n
=\bigcup_{l=0}^{\lfloor1/\sqrt{\delta}\rfloor+1} D_{r+\delta} \left(\sigma_0+i \left(n+l\sqrt{\delta}\right)\right).
\end{gather*}
We take $s_n=\sigma_n+i t_n \in E_n$ so that $|f(z)|$ takes the maximum value at $s_n$. 
The function $F(x)=\log^2(x+3)$ is convex for $x \geq0$. 
Hence we have 
\begin{align}\label{e6.9}
\sum_{n=\lfloor T\rfloor}^{\lfloor2T\rfloor+1}\log^2\left(2 \max_{z \in E_n}|f(z)|+3\right)
&=\sum_{n=\lfloor T\rfloor}^{\lfloor2T\rfloor+1} F(2|f(s_n)|)\\
&\ll T F\left(\frac{1}{T} \sum_{n=\lfloor T\rfloor}^{\lfloor2T\rfloor+1} 2|f(s_n)|\right) \nonumber\\
&\leq T F\left(\frac{2}{T^{\frac{1}{2}}} \left(\sum_{n=\lfloor T\rfloor}^{\lfloor2T\rfloor+1} 
|f(s_n)|^2\right)^{1/2}\right). \nonumber
\end{align}
The remaining work is to estimate $\sum_n|f(s_n)|^2$. 
For this, we define  
\begin{gather*}
S_j
=\left\{ s_n \mid n \equiv j \bmod \left(4\lfloor r+2\delta \rfloor+6\right)\right\}.
\end{gather*}
Then we have 
\begin{equation}\label{e6.10}
\sum_{n=\lfloor T\rfloor}^{\lfloor2T\rfloor+1} |f(s_n)|^2
\ll_{\sigma_1} \sum_{\substack{s_n \in S_j \\ \lfloor T\rfloor \leq n \leq \lfloor2T\rfloor+1}} |f(s_n)|^2.
\end{equation}
Note that $f(z)$ is regular for $|z-(\sigma_0+it_n)| \leq r+2\delta$ and $|s_n-(\sigma_0+it_n)| \leq r+\delta$.
Hence, the inequality of \cite[Lemma in p.\,256]{Titchmarsh1986} deduces 
\begin{gather*}
|f(s_n)|^2
\leq \frac{1}{\pi\delta^2} \iint_{D_{r+2\delta}(\sigma_0+it_n)} |f(z)|^2 \,dxdy.
\end{gather*}
Moreover, if $s_m, s_n \in S_j$ with $m>n$, then we have 
\begin{gather*}
|t_m-t_n|
\geq \{m-(r+\delta)\}-\{n+(\lfloor \tfrac{1}{\sqrt{\delta}} \rfloor+1) \sqrt{\delta} 
+(r+\delta)\}
>2r + 4\delta.
\end{gather*}
Thus $D_{r+2\delta}(\sigma_0+it_m)$ and $D_{r+2\delta}(\sigma_0+it_m)$ are disjoint, and therefore,  
\begin{align}\label{e6.11}
\sum_{\substack{s_n \in S_j \\ \lfloor T \rfloor \leq n \leq \lfloor 2T \rfloor+1}} |f(s_n)|^2
&\leq \frac{1}{\pi \delta^2} \int_{\sigma_0-(r+2\delta)}^{\sigma_0+(r+2\delta)}
\int_{T-1-(2r+3\delta)}^{2T+3+(2r+3\delta)} |f(x+iy)|^2 \,dydx \\
&\ll_{\alpha,\sigma_1,\epsilon_1} \int_{\frac{1}{2}+c_1 \epsilon_1}^{2\sigma_0+2}
\int_{T-1-(2r+3\delta)}^{2T+3+(2r+3\delta)} |L(\lambda,\alpha,x+iy)|^2 \,dydx \nonumber
\end{align}
for a positive constant $c_1=c_1(\sigma_1)$. 
With regard to the inner integral, we recall that for all $\sigma\geq1/2+\epsilon_1$, 
\begin{gather*}
\int_{0}^{T} |L(\lambda,\alpha,\sigma+it)|^2 \,dt
\ll_{\lambda,\alpha,\epsilon_1} T
\end{gather*}
for $T \geq T_0$ with a positive constant $T_0=T_0(\lambda,\alpha,\epsilon_1)$ by \cite[Theorem 3.3.1]{LaurincikasGarunkstis2002}. 
Hence we see that 
\begin{equation}\label{e6.12}
\int_{\frac{1}{2}+c_1 \epsilon_1}^{2\sigma_0+2} \int_{T-1-(2r+3\delta)}^{2T+3+(2r+3\delta)} 
|L(\lambda,\alpha,x+iy)|^2 \,dydx
\ll_{\lambda,\alpha,\sigma_1,\epsilon_1} T.
\end{equation}
From the above \eqref{e6.7}--\eqref{e6.12}, we conclude 
\begin{gather*}
\int_{T}^{2T} \left(\sum_{|\rho-s_0(t)| \leq r_1} \log|\sigma+it-\rho|\right)^2 \,dt 
\ll_{\lambda,\alpha,\sigma_1,\epsilon_1} T
\end{gather*}
as desired. 
\end{proof}

For the second integral of the right-hand side of \eqref{e6.5}, we obtain the following estimate. 

\begin{lemma}\label{l6.5}
Let $\lambda \in \mathbb{R}$ and $0<\alpha \leq1$. 
Let $\sigma_1$ be a large fixed positive constant. 
Let $\epsilon_1>0$ be a small fixed real number. 
Then there exists $T_0=T_0(\lambda,\alpha, \sigma_1, \epsilon_1)>0$ such that for all $T \geq T_0$ and for all $1/2+\epsilon_1 \leq \sigma \leq \sigma_1$, we have 
\begin{gather*}
\int_{T}^{2T} \log^2 M_r(t) \,dt 
\ll_{\lambda,\alpha,\sigma_1,\epsilon_1} T.
\end{gather*}
\end{lemma}

\begin{proof}
We find that the function $G(x)=\log^2{x}$ is convex for $x>e$ and $M_r(t)^2>e$. 
Then, applying Jensen's inequality, we have 
\begin{gather*}
\int_{T}^{2T} \log^2 M_r(t) \,dt 
\ll \int_{T}^{2T} \log^2 M_r(t)^2 \,dt
\leq T \log^2\left(\frac{1}{T} \int_{T}^{2T} M_r(t)^2 \,dt \right). 
\end{gather*}
Also we have 
\begin{gather*}
\int_{T}^{2T} M_r(t)^2 \,dt
\leq \sum_{n=\lfloor T \rfloor}^{\lfloor 2T \rfloor+1} \int_{n}^{n+1} M_r(t)^2 \,dt.
\end{gather*}
Let 
\begin{gather*}
F_n
=\bigcup_{n \leq t \leq n+1} D_r(s_0(t)),
\end{gather*}
and let $s'_n \in F_n$ be a point at which $|f(z)|$ takes the maximum value. 
Then we have 
\begin{gather*}
\sum_{n=\lfloor T \rfloor}^{\lfloor 2T \rfloor+1} \int_{n}^{n+1} M_r(t)^2 \,dt
\ll \sum_{n=\lfloor T \rfloor}^{\lfloor 2T \rfloor+1} |f(s'_n)|^2
+T,
\end{gather*}
and by a similar argument in the proof of Lemma \ref{l6.3}, we see that 
\begin{gather*}
\sum_{n=\lfloor T \rfloor}^{\lfloor 2T \rfloor+1} |f(s'_n)|^2
\ll_{\lambda,\alpha,\sigma_1,\epsilon_1} T
\end{gather*}
for $T \geq T_0$ with a positive constant $T_0=T_0(\lambda,\alpha,\epsilon_1)$. 
Hence the result follows. 
\end{proof}

By Lemmas \ref{l6.4} and \ref{l6.5}, we have 
\begin{gather*}
\int_{T}^{2T} \log^2|L(\lambda,\alpha,\sigma+it)| \,dt 
\ll T
\end{gather*}
for large $T$. 
To finish the proof of Proposition \ref{p6.2}, we need an estimate for small $T$. 

\begin{lemma}\label{l6.6}
Let $\lambda \in \mathbb{R}$ and $0<\alpha \leq1$. 
Let $\sigma_1$ be a large fixed positive constant. 
Let $\epsilon_1>0$ be a small fixed real number. 
Then for any fixed constant $T_0>0$, and for all $1/2+\epsilon_1 \leq \sigma \leq \sigma_1$, we have 
\begin{equation}\label{e6.13}
\int_{0}^{T_0} \log^2|L(\lambda,\alpha,\sigma+it)| \,dt 
\ll 1.
\end{equation}
The implied constant depends only on $\lambda$, $\alpha$, $\sigma_1$, $\epsilon_1$, and $T_0$.
\end{lemma}

\begin{proof}
Let $\rho_1,\ldots,\rho_n$ be all zeros of $L(\lambda,\alpha,s)$ in the rectangle $1/2+\epsilon_1 \leq \sigma \leq \sigma_1$, $0 \leq t \leq T_0$. 
Here $n \geq 0$ is a finite integer determined from $\lambda$, $\alpha$, $\sigma_1$, $\epsilon_1$, and $T_0$. 
We take a positive real number $r_0=r_0(\lambda,\alpha,\sigma_1,\epsilon_1,T_0)$ so that the disks $D_{r_0}(\rho_1),\ldots, D_{r_0}(\rho_n)$ are distinct. 
Let $\Omega$ be the closure of the set
\begin{gather*}
\{1/2+\epsilon_1 \leq \sigma \leq \sigma_1,~ 0 \leq t \leq T_0\} 
\setminus (D_{r_0}(\rho_1) \cup \cdots \cup D_{r_0}(\rho_n))
\end{gather*}
and define 
\begin{gather*}
M
=M(\lambda,\alpha,\sigma_1,\epsilon_1,T_0)
=\max_{s \in \Omega} \log^2|L(\lambda,\alpha,s)|.
\end{gather*}
We divide the integral in \eqref{e6.13} into 
\begin{gather*}
\int_{I_1} \log^2|L(\lambda,\alpha,\sigma+it)| \,dt
\quad\text{and}\quad
\int_{I_2} \log^2|L(\lambda,\alpha,\sigma+it)| \,dt,
\end{gather*}
where 
\begin{align*}
I_1
&=\{\sigma+it \mid 0 \leq t \leq T_0\} \cap \Omega, \\
I_2
&=\{\sigma+it \mid 0 \leq t \leq T_0\} \cap \Omega^c.
\end{align*}
The integral over $I_1$ is estimated as 
\begin{equation}\label{e6.14}
\int_{I_1} \log^2|L(\lambda,\alpha,\sigma+it)| \,dt
\leq M T_0
\ll1,
\end{equation}
where the implied constant depends only on $\lambda$, $\alpha$, $\sigma_1$, $\epsilon_1$, and $T_0$. 
Let 
\begin{gather*}
L(\lambda,\alpha,s)
=(s-\rho_k)^{m_k} L_1(\lambda,\alpha,s),
\end{gather*}
where $m_k$ is the order of the zero at $s=\rho_k$.
Then we have 
\begin{gather*}
\log^2|L(\lambda,\alpha,s)|
\ll m_k^2 \log^2|s-\rho_k|
+\log^2|L_1(\lambda,\alpha,s)|.
\end{gather*}
Therefore, the integral over $I_2$ is 
\begin{align}\label{e6.15}
&\int_{I_2} \log^2|L(\lambda,\alpha,\sigma+it)| \,dt\\
&\leq \sum_{k=1}^{n} m_k^2 \int_{\gamma_k-r_0}^{\gamma_k+r_0} \log^2|\sigma+it-\rho_k| \,dt
+\sum_{k=1}^{n} \int_{\gamma_k-r_0}^{\gamma_k+r_0} \log^2|L_1(\lambda,\alpha,\sigma+it)| \,dt, \nonumber
\end{align}
where $\rho_k=\beta_k+i\gamma_k$.
We see that 
\begin{gather*}
\int_{\gamma_k-r_0}^{\gamma_k+r_0} \log^2|\sigma+it-\rho_k| \,dt
\leq 2 \int_{0}^{r_0} \log^2(x+|\sigma-\beta_k|) \,dt
\ll 1
\end{gather*}
and 
\begin{gather*}
\int_{\gamma_k-r_0}^{\gamma_k+r_0} 
\log^2|L_1(\lambda,\alpha,\sigma+it)| \,dt 
\ll1
\end{gather*}
with the implied constant depending on $\lambda$, $\alpha$, $\sigma_1$, $\epsilon_1$, and $T_0$. 
The maximum value of $m_k$ is determined only from $\lambda$, $\alpha$, $\sigma_1$, $\epsilon_1,T_0$, and hence we have 
\begin{equation}\label{e6.16}
\int_{I_2} \log^2|L(\lambda,\alpha,\sigma+it)| \,dt
\ll1
\end{equation}
due to \eqref{e6.15}. 
Lemma \ref{l6.6} is deduced from estimates \eqref{e6.14} and \eqref{e6.16}. 
\end{proof}

\begin{proof}[Proof of Proposition \ref{p6.2}]
By \eqref{e6.5} and Lemmas \ref{l6.4} and \ref{l6.5}, we have 
\begin{gather*}
\int_{2^{k-1}T_0}^{2^kT_0} \log^2|L(\lambda,\alpha,\sigma+it)| \,dt 
\ll 2^{k-1} T_0
\end{gather*}
for $k \geq1$.
Together with Lemma \ref{l6.6}, by summing up over $k$, we have the result. 
\end{proof}

\subsection{Proof of Proposition \ref{p6.1}}\label{s6.2}
We begin with constructing certain functions that approximate $\log|z|$ by following the way in \cite{Guo1996b}. 
Let 
\begin{gather*}
f(u)=
\begin{cases}
\displaystyle{\exp\left( -(b-a) \left(\frac{1}{u-a}+\frac{1}{b-u}\right) \right)}
&\quad\text{if $a<u<b$}, \\
0
&\quad\text{otherwise}
\end{cases}
\end{gather*}
with $a, b \in \mathbb{R}$ and $0<b-a<1$. 
Then we define
\begin{gather*}
F(x)
=\frac{\displaystyle \int_{-\infty}^{x+b-a_1} f(u) \,du}{\displaystyle \int_{-\infty}^{\infty} f(u) \,du} 
\cdot \frac{\displaystyle \int_{x+a-b_1}^\infty f(u) \,du}{\displaystyle \int_{-\infty}^{\infty} f(u) \,du}
\end{gather*}
for $a_1, b_1 \in \mathbb{R}$ with $a_1<b_1$. 
With the above setting, we define
\begin{gather*}
\Phi(z)
=F(|z|) \log|z|,
\end{gather*}
which is infinitely differentiable and supported on $a_2 \leq |z| \leq b_2$ with $a_2=a_1-(b-a)$ and $b_2=b_1+(b-a)$. 
Moreover, we take the above $a$, $b$, $a_1$, $b_1$, $a_2$, $b_2$ as the functions 
\begin{gather*}
\begin{array}{ll}
a=1-(\log{T})^{-\gamma}, & b=1, \\
a_1=2(\log{T})^{-\gamma}, & b_1=(\log{T})^\gamma, \\
a_2=(\log{T})^{-\gamma}, & b_2=(\log{T})^\gamma+(\log{T})^{-\gamma}
\end{array}
\end{gather*}
with a constant $0<\gamma<1$. 
Then the following lemma holds.  

\begin{lemma}[Lemma 3.1.1 of \cite{Guo1996b}]\label{l6.7}
The function $\Phi(z)$ satisfies the following conditions. 
\begin{enumerate}
\item[$\mathrm{(1)}$] 
$\Phi(z)=\log|z|$ for $a_1 \leq|z| \leq b_1$.  
\item[$\mathrm{(2)}$] 
$|\Phi(z)| \leq |\log|z||$ for $a_2 \leq |z| \leq a_1$, $b_1 \leq |z| \leq b_2$. 
\item[$\mathrm{(3)}$] 
The Fourier transform of $\Phi(z)$ is estimated as 
\begin{gather*}
\widetilde{\Phi}(z)
\ll \left(|\log a_2|+|\log b_2|\right) 
\min \left(b_2^2, \frac{b_2^2+a_2^{-2}}{(b-a)^2 x^2}, 
\frac{b_2^2+a_2^{-2}}{(b-a)^2 y^2}, \frac{b_2^2+a_2^{-2}}{(b-a)^4 x^2 y^2} \right)
\end{gather*}
for large $T$ with $z=x+iy$. 
\end{enumerate}
\end{lemma}

\begin{proof}[Proof of Proposition \ref{p6.1}]
Since the function $\Phi(z)$ belongs to the class $\Lambda$, we have 
\begin{equation}\label{e6.17}
\frac{1}{T} \int_{0}^{T} \Phi(L(\lambda, \alpha, \sigma+it)) \,dt
=\int_{\mathbb{C}} \Phi(z) M_\sigma(z; \alpha) \,|dz|
+E,
\end{equation}
by Theorem \ref{t3.2}, where 
\begin{gather*}
E
\ll_{\lambda, \sigma_1} \exp\left(-\frac{1}{2} (\log{T})^{\theta/2}\right) 
\int_{\Omega} \left|\widetilde{\Phi}(z)\right| \,|dz|
+\int_{\mathbb{C} \setminus \Omega} \left|\widetilde{\Phi}(z)\right| \,|dz|.
\end{gather*}
By the first inequality of condition $\mathrm{(3)}$ in Lemma \ref{l6.7}, the first integral is estimated as 
\begin{gather*}
\int_{\Omega} \left|\widetilde{\Phi}(z)\right| \,|dz|
\ll \int_{-(\log{T})^\delta}^{(\log{T})^\delta} \int_{-(\log{T})^\delta}^{(\log{T})^\delta} (\log{T})^{2\gamma+\epsilon} \,dxdy
\ll (\log{T})^{2\gamma+2\delta+\epsilon}
\end{gather*}
with sufficiently small $\epsilon>0$. 
For the estimate of the second integral, we consider a covering of the region $\mathbb{C}\backslash\Omega$ as follows: 
\begin{gather*}
\mathbb{C} \setminus \Omega 
\subset U_1 \cup U_2 \cup U_3, 
\end{gather*}
where 
\begin{align*}
U_1
&=\left\{ z=x+iy ~\middle|~ |x| \geq (\log{T})^\delta,~ |y| \leq (\log{T})^{\delta/2} \right\}, \\
U_2
&=\left\{ z=x+iy ~\middle|~ |x| \leq (\log{T})^{\delta/2},~ |y| \geq (\log{T})^\delta \right\}, \\
U_3
&=\left\{ z=x+iy ~\middle|~ |x| \geq (\log{T})^{\delta/2},~ |y| \geq (\log{T})^{\delta/2} \right\}. 
\end{align*}
Then we have 
\begin{gather*}
\int_{U_1} \left|\widetilde{\Phi}(z)\right| \,|dz|
\ll \int_{-(\log{T})^{\delta/2}}^{(\log{T})^{\delta/2}} \int_{(\log{T})^\delta}^{\infty} 
\frac{(\log{T})^{4\gamma+\epsilon}}{x^2} \,dxdy
\ll (\log{T})^{4\gamma-{\delta/2}+\epsilon}
\end{gather*}
by the second inequality of condition $\mathrm{(3)}$. 
Similarly we have  
\begin{gather*}
\int_{U_2} \left|\widetilde{\Phi}(z)\right| \,|dz|
\ll (\log{T})^{4\gamma-{\delta/2}+\epsilon}
\end{gather*}
by using the third inequality of condition $\mathrm{(3)}$, and furthermore, the fourth inequality deduces 
\begin{gather*}
\int_{U_3} \left|\widetilde{\Phi}(z)\right| \,|dz|
\ll \int_{(\log{T})^{\delta/2}}^{\infty} \int_{(\log{T})^{\delta/2}}^{\infty} 
\frac{(\log{T})^{6\gamma+\epsilon}}{x^2y^2} \,dxdy
\ll (\log{T})^{6\gamma-\delta+\epsilon}. 
\end{gather*}
Assuming $\gamma<\delta/4$, we have $4\gamma-\delta/2>6\gamma-\delta$. 
Thus we obtain 
\begin{align}\label{e6.18}
E
&\ll_{\lambda,\sigma_1} \exp\left(-\frac{1}{2} (\log{T})^{\theta/2}\right) (\log{T})^{2\gamma+2\delta+\epsilon}
+(\log{T})^{4\gamma-\delta/2+\epsilon}\\
&\ll(\log{T})^{4\gamma-\delta/2+\epsilon}\nonumber
\end{align}
for any $T \geq T_0$ with a positive constant $T_0=T_0(\alpha,\theta,\delta,\sigma_1,\gamma)$. 

Next, we consider the difference
\begin{gather*}
E_L
=\frac{1}{T} \int_{0}^{T} \Phi(L(\lambda, \alpha, \sigma+it)) \,dt
-\frac{1}{T} \int_{0}^{T} \log|L(\lambda, \alpha, \sigma+it)| \,dt
\end{gather*}
coming from the left-hand side of \eqref{e6.17}. 
Let 
\begin{align*}
I
&=\{ t \in [0,T] \mid |L(\lambda, \alpha, \sigma+it)| \leq a_1 \}, \\
J
&=\{ t \in [0,T] \mid |L(\lambda, \alpha, \sigma+it)| \geq b_1 \}. 
\end{align*}
Then Lemma \ref{l6.7} gives 
\begin{gather*}
|E_L| 
\leq \frac{1}{T} \int_{I \cup J} |\log|L(\lambda, \alpha, \sigma+it)|| \,dt.
\end{gather*}
Applying Cauchy's inequality, we see that this is
\begin{equation}\label{e6.19}
\leq \left(\frac{\mu_1(I)+\mu_1(J)}{T}\right)^{1/2} 
\left(\frac{1}{T} \int_{0}^{T} \log|L(\lambda, \alpha, \sigma+it)|^2 \,dt\right)^{1/2}.
\end{equation}
Let
\begin{align*}
R
&=\{ z=x+iy \mid |x| \leq (\log{T})^{-\gamma},~ |y| \leq (\log|{T})^{-\gamma}\}, \\
R'
&=\{ z=x+iy \mid |x| \leq \tfrac{1}{\sqrt{2}} (\log{T})^\gamma,~ |y| \leq \tfrac{1}{\sqrt{2}} (\log{T})^\gamma\}.
\end{align*}
Since we have $(\log{T})^{-\gamma} \gg (\log{T})^{-1/16}$, it is deduced from Theorem \ref{t1.7} that 
\begin{align}\label{e6.20}
\frac{\mu_1(I)}{T}
&\ll_{\lambda,\alpha,\sigma_1,\epsilon_1} \int_{R} M_\sigma(z;\alpha) \,|dz|
+(\mu_2(R)+1) (\log{T})^{-1/8} \\
&\ll \int_{R} M_\sigma(z;\alpha) \,|dz|
+(\log{T})^{-1/8} \nonumber
\end{align}
and
\begin{align}\label{e6.21}
\frac{\mu_1(J)}{T}
&=1-\frac{\mu_1([0,T] \setminus J)}{T} \\
&\ll_{\lambda,\alpha,\sigma_1,\epsilon_1} 1-\int_{R'} M_\sigma(z;\alpha) \,|dz|
+(\mu_2(R')+1) (\log{T})^{-1/8} \nonumber\\
&\ll \int_{\mathbb{C} \setminus R'} M_\sigma(z;\alpha) \,|dz|
+(\log{T})^{2\gamma-1/8}.\nonumber
\end{align}
for any $T \geq T^{(1)}_0$ with a positive constant $T^{(1)}_0=T^{(1)}_0(\alpha,\sigma_1,\epsilon_1)$. 
The integrals in \eqref{e6.20} and \eqref{e6.21} are estimated as 
\begin{equation}\label{e6.22}
\int_{R} M_\sigma(z;\alpha) \,|dz|
\ll_{\alpha,\epsilon_1} \int_{R} 1 \,dxdy
\ll (\log{T})^{-2\gamma}
\end{equation}
and
\begin{equation}\label{e6.23}
\int_{\mathbb{C} \setminus R'} M_\sigma(z;\alpha) \,|dz|
\ll_{\alpha,\epsilon_1} \int_{\tfrac{1}{\sqrt{2}} (\log{T})^\gamma}^{\infty} e^{-r^2} r \,dr
\ll \exp\left(-\frac{1}{4} (\log{T})^{2\gamma}\right)
\end{equation}
due to Proposition \ref{p2.5}. 
We then use Proposition \ref{p6.2}. 
It follows from Proposition \ref{p6.2} and estimates \eqref{e6.19}--\eqref{e6.23} that $E_L$ is estimated as 
\begin{equation}\label{e6.24}
E_L
\ll_{\lambda,\alpha,\sigma_1,\epsilon_1} (\log{T})^{2\gamma-1/8}
\end{equation}
for any $T \geq T^{(2)}_0$ with a constant $T^{(2)}_0=T^{(2)}_0(\alpha,\sigma_1,\epsilon_1,\gamma) \geq T^{(1)}_0$. 

Finally, we estimate 
\begin{gather*}
E_R
=\int_{\mathbb{C}} \Phi(z) M_\sigma(z; \alpha) \,|dz|
-\int_{\mathbb{C}} \log|z| M_\sigma(z; \alpha) \,|dz|.
\end{gather*}
It is estimated as 
\begin{gather*}
|E_R|
\leq \int_{|z| \leq a_1} |\log|z| |M_\sigma(z; \alpha) \,|dz|
+\int_{|z| \geq b_1} |\log|z| |M_\sigma(z; \alpha) \,|dz|
\end{gather*}
by Lemma \ref{l6.7}. 
Again applying Proposition \ref{p2.5}, we see that the first integral is 
\begin{gather*}
\ll_{\alpha,\epsilon_1} |\log a_1| a_1^2 
\ll (\log{T})^{-2\gamma+\epsilon}, 
\end{gather*}
and the second integral is 
\begin{gather*}
\ll_{\alpha,\epsilon_1} \int_{b_1}^{\infty} \log{r} e^{-r^{2}} r \,dr
\ll \int_{b_1}^{\infty} r^{-2} \,dr
\ll (\log{T})^{-\gamma}. 
\end{gather*}
If we take $\epsilon<\gamma$, we have 
\begin{equation}\label{e6.25}
E_R
\ll_{\alpha,\epsilon_1} (\log{T})^{-\gamma}. 
\end{equation}
We take $\gamma=\delta/12$ and fix $\theta, \delta>0$ with $\theta+\delta<1/4$. 
Combining estimates \eqref{e6.18}, \eqref{e6.24}, and \eqref{e6.25}, we conclude that 
\begin{gather*}
\frac{1}{T} \int_{0}^{T} \log|L(\lambda, \alpha, \sigma+it)| \,dt
-\int_{\mathbb{C}} \log|z| M_\sigma(z; \alpha) \,|dz|
\ll_{\lambda,\alpha,\sigma_1,\epsilon_1} (\log{T})^{-\frac{\delta}{12}}.
\end{gather*}
Hence the result follows. 
\end{proof}

\subsection{Completion of the proof}\label{s6.3}
We use the following lemma to prove Theorem \ref{t1.8}. 

\begin{lemma}[Lemma 8.4.11 of \cite{LaurincikasGarunkstis2002}]\label{l6.8}
Let $1/2 \leq \sigma \leq 1+\alpha$. 
Then we have  
\begin{gather*}
2\pi \int_{\sigma}^{1+\alpha} N(T, u; \lambda, \alpha) \,du
=\sigma T \log{\alpha}
+\int_{0}^{T} \log|L(\lambda, \alpha, \sigma+it)| \,dt 
+O(\log{T})
\end{gather*}
for sufficiently large $T$, where $N(T, u; \lambda, \alpha)$ is the number of zeros of $L(\lambda, \alpha, \sigma+it)$ in the region $\sigma>u$, $0<t<T$. 
\end{lemma}

\begin{proof}[Proof of Theorem \ref{t1.8}]
Let $\sigma>1/2$ be a fixed real number. 
By Proposition \ref{p6.1} and Lemma \ref{l6.8}, we have  
\begin{align*}
&2\pi \int_{\sigma}^{1+\alpha} N(T, u; \lambda, \alpha) \,du \\
&=\sigma T \log{\alpha}
+T \int_{\mathbb{C}} \log|z| M_\sigma(z; \alpha) \,|dz|
+O_{\lambda,\alpha,\sigma}\left(T (\log{T})^{-A}\right)
\end{align*}
for $T \geq T_0$. 
If a small positive real number $h$ satisfies 
\begin{gather*}
\frac{1}{2}
<\frac{1}{2}\left(\sigma+\frac{1}{2}\right)
\leq \sigma-h
<\sigma
<\sigma+h
\leq\sigma+1, 
\end{gather*}
we have also
\begin{align*}
&2\pi \int_{\sigma\pm h}^{1+\alpha} N(T, u; \lambda, \alpha) \,du \\
&=\sigma T \log{\alpha} 
+T \int_{\mathbb{C}} \log|z| M_{\sigma\pm h}(z; \alpha) \,|dz|
+O_{\lambda,\alpha,\sigma}\left(T (\log{T})^{-A}\right).
\end{align*}
Hence the formula
\begin{align}\label{e6.26}
&2\pi \int_{\sigma}^{\sigma+h} N(T, u; \lambda, \alpha)du \\
&=h T \log{\alpha}
+T (\phi_\alpha(\sigma+h)-\phi_\alpha(\sigma))
+O_{\lambda,\alpha,\sigma} \left(T (\log{T})^{-A}\right) \nonumber 
\end{align}
holds, where
\begin{gather*}
\phi_\alpha(\sigma)
=\int_{\mathbb{C}} \log|z| M_\sigma(z; \alpha) \,|dz|. 
\end{gather*}
By Corollary \ref{c5.4}, we see that $\phi_\alpha(\sigma)$ is infinitely differentiable, and therefore 
\begin{gather*}
\frac{\partial^n}{\partial \sigma^n} \phi_\alpha(\sigma)
=\int_{\mathbb{C}} \log|z| \frac{\partial^n}{\partial \sigma^n} M_\sigma(z; \alpha) \,|dz|.
\end{gather*}
Hence we have 
\begin{gather*}
\phi_\alpha(\sigma+h)-\phi_\alpha(\sigma)
=h \int_{\mathbb{C}} \log|z| \frac{\partial}{\partial \sigma} M_\sigma(z; \alpha) \,|dz|
+O_{\alpha,\sigma}(h^2).
\end{gather*}
Since the function $N(T, u; \lambda, \alpha)$ is decreasing in $u$, we have, by \eqref{e6.26},
\begin{align*}
&N(T, \sigma; \lambda, \alpha) \\
&\geq T \frac{\log{\alpha}}{2\pi}
+\frac{T}{2\pi} \int_{\mathbb{C}} \log|z| \frac{\partial}{\partial \sigma} M_\sigma(z; \alpha) \,|dz|
+O_{\lambda,\alpha,\sigma} \left( h^{-1} T (\log{T})^{-A} + h T \right).
\end{align*}
Similarly, we have 
\begin{align*}
&N(T, \sigma; \lambda, \alpha) \\
&\leq T \frac{\log{\alpha}}{2\pi}
+\frac{T}{2\pi} \int_{\mathbb{C}} \log|z| \frac{\partial}{\partial \sigma} M_\sigma(z; \alpha) \,|dz|
+O_{\lambda,\alpha,\sigma}\left( h^{-1} T (\log{T})^{-A} + h T \right)
\end{align*}
by considering $\sigma-h$ instead of $\sigma+h$.
We take $h=(\log{T})^{-A/2}$. 
The above error terms are $\ll T (\log{T})^{-A/2}$, and hence we have 
\begin{gather*}
N(T, \sigma; \lambda, \alpha)
=T \frac{\log{\alpha}}{2\pi}
+\frac{T}{2\pi} \int_{\mathbb{C}} \log|z| \frac{\partial}{\partial \sigma} M_\sigma(z; \alpha) \,|dz|
+O_{\lambda,\alpha,\sigma}\left(T (\log T)^{-A/2}\right). 
\end{gather*}
Thus we obtain for any fixed $1/2<\sigma_1<\sigma_2$,
\begin{align*}
&N(T, \sigma_1,\sigma_2; \lambda, \alpha) \\
&=N(T, \sigma_2; \lambda, \alpha)-N(T, \sigma_1; \lambda, \alpha)\\
&=\frac{T}{2\pi} \int_{\mathbb{C}} \log|z| 
\left(\frac{\partial}{\partial \sigma} M_{\sigma_2}(z; \alpha)
-\frac{\partial}{\partial \sigma} M_{\sigma_1}(z; \alpha)\right) \,|dz|
+O\left(T (\log{T})^{-A/2}\right)\\
&=\frac{T}{2\pi} \int_{\sigma_1}^{\sigma_2} \int_{\mathbb{C}} 
\log|z| \frac{\partial^2}{\partial \sigma^2} M_\sigma(z; \alpha) \,|dz| \,d\sigma
+O\left(T (\log{T})^{-A/2}\right)
\end{align*}
as desired, where the implied constant depends only on $\lambda$, $\alpha$, $\sigma_1$ and $\sigma_2$. 
\end{proof}

\appendix
\section{The study of Borchsenius and Jessen}\label{sa}
In this paper, we have proved an asymptotic formula for zeros of $L(\lambda,\alpha,s)$ of the form 
\begin{gather*}
N(T, \sigma_1, \sigma_2; \lambda, \alpha)
=C T
+O\left(T (\log{T})^{-A}\right),
\qquad 
T \to\infty
\end{gather*} 
in the case of $\alpha \in \mathfrak{S}$. 
As we have remarked in Section \ref{s1.2}, we may also obtain the formula
\begin{gather*}
C(\sigma_1, \sigma_2;\alpha)
=\lim_{T \to\infty} \frac{1}{T} N(T, \sigma_1, \sigma_2; \lambda, \alpha)
\end{gather*}
if $\alpha$ satisfies condition \eqref{c1} of Definition \ref{d1.6}. 
Since the proof is a simple analogue of the method of Borchsenius and Jessen \cite{BorchseniusJessen1948}, we follow it in this section as an appendix of this paper. 
For the proof, we refer to three results from \cite{BorchseniusJessen1948} as follows. 

\begin{lemma}\label{la.1}
Let $-\infty \leq \alpha<\alpha_0<\beta_0<\beta \leq +\infty$ and $-\infty<\gamma_0<+\infty$, and let $f_1(s), f_2(s), \ldots$ be a sequence of functions almost periodic in $[\alpha, \beta]$ converging uniformly in $[\alpha_0, \beta_0]$ towards a function $f(s)$. 
Suppose that none of the functions is identically zero. 
Suppose further, that $f(s)$ is continued as a regular function in the half strip $\alpha<\sigma<\beta$, $t>\gamma_0$, and that for any fixed $\alpha<\alpha_1<\beta_1<\beta$ and $\gamma>\gamma_0$, 
\begin{gather*}
\limsup_{\delta \to\infty} \frac{1}{\delta-\gamma} \int_{\gamma}^{\delta}
\left\{ \int_{\alpha_1}^{\beta_1} |f(\sigma+it)-f_N(\sigma+it)|^p \,d\sigma \right\}^{1/p} \,dt 
\to0 
\end{gather*}
as $N \to\infty$ with an index $p>0$. 

Then the function 
\begin{gather*}
\phi_f(\sigma;\gamma,\delta)
=\frac{1}{\delta-\gamma} \int_{\gamma}^{\delta} \log|f(\sigma+it)| \,dt
\end{gather*}
converges as $\delta \to\infty$ for any fixed $\gamma>\gamma_0$ uniformly in $[\alpha, \beta]$ towards a function $\phi_f(\sigma)$, which is called the Jensen function of $f(s)$. 
The Jensen function $\phi_{f_N}(\sigma)$ of $f_N(s)$ converges as $N \to\infty$ uniformly in $[\alpha, \beta]$ towards $\phi_f(\sigma)$. 

Moreover, for every strip $(\sigma_1, \sigma_2)$ with $\alpha<\sigma_1<\sigma_2<\beta$, if $\phi_f(\sigma)$ is differentiable at $\sigma_1$ and $\sigma_2$, then the limit value 
\begin{equation}\label{ea.1}
C_f(\sigma_1,\sigma_2)
=\lim_{\delta \to\infty} \frac{N_f(\sigma_1,\sigma_2;\gamma,\delta)}{\delta-\gamma}
\end{equation}
exists for any fixed $\gamma>\gamma_0$ and determined by the formula
\begin{gather*}
C_f(\sigma_1,\sigma_2)
=\frac{1}{2\pi} (\phi'_f(\sigma_2)-\phi'_f(\sigma_1)),
\end{gather*}
where $N_f(\sigma_1,\sigma_2;\gamma,\delta)$ denotes the number of zeros of $f(s)$ in the rectangle $\sigma_1<\sigma<\sigma_2$, $\gamma<t<\delta$. 
\end{lemma}

\begin{proof}
The first part of the conclusion is from \cite[Theorem 1]{BorchseniusJessen1948}, and the second part is seen in the argument following the statement of \cite[Theorem 1]{BorchseniusJessen1948}. 
\end{proof}

\begin{lemma}[Theorem 5 of \cite{BorchseniusJessen1948}]\label{la.2}
Let $\ell(z)=\ell_1z+\ell_2z^2+\cdots$ and $m(z)=m_1z+m_2z^2+\cdots$ be power series converging in a circle $|z|<\rho$, and such that $\ell_1 \neq0$ and $m_1 \neq0$. 
Let $r_0,r_1,\ldots$ be a sequence of real numbers such that $0<r_n<\rho$ for all $n$, and let $\lambda_0,\lambda_1,\ldots$ be a sequence of real numbers differing from each other and from zero. 

For every $N$, we define
\begin{gather*}
f_N(t_0,\ldots,t_N)
=\sum_{n=0}^{N} \ell(r_nt_n)
\quad\text{and}\quad
g_N(t_0,\ldots,t_N)
=\sum_{n=0}^{N} \lambda_n m(r_nt_n),
\end{gather*}
where $(t_0,\ldots,t_N)$ describes the torus $T^{N+1}=\prod_{n=0}^{N} \mathbb{C}^1$, where $\mathbb{C}^1=\{z \in \mathbb{C} \mid |z|=1\}$.

Let $m_N$ be the normalized Haar measure of $T^{N+1}$. 
We define a probability measure $\mu_N$ on $(\mathbb{C},\mathcal{B}(\mathbb{C}))$ as 
\begin{gather*}
\mu_N(A)
=\int_{T^{N+1}} 1_A(f_N(t_0,\ldots,t_N)) \,dm_N(t_0,\ldots,t_N). 
\end{gather*}
We also define $\nu_N$ as a measure on $(\mathbb{C},\mathcal{B}(\mathbb{C}))$ such that
\begin{gather*}
\nu_N(A)
=\int_{T^{N+1}} 1_A(f_N(t_0,\ldots,t_N)) |g_N(t_0,\ldots,t_N)|^2 \,dm_N(t_0,\ldots,t_N). 
\end{gather*}

Then, if $r_n \to 0$, the measures $\mu_N$ and $\nu_N$ are absolutely continuous with continuous density functions $F_N(z)$ and $G_N(z)$, respectively, for $N \geq N_0$ with a positive integer $N_0$, and $F_N(z)$ and $G_N(z)$ possess continuous partial derivatives of order $\leq p$ for $N \geq N_p$ with a positive integer $N_p$. 

Moreover, if the series 
\begin{gather*}
S_0
=\sum_{n=0}^{\infty} r_n^2, 
\quad
S_1
=\sum_{n=0}^{\infty} |\lambda_n| r_n^2, 
\quad
S_2
=\sum_{n=0}^{\infty} \lambda_n^2 r_n^2
\end{gather*}
converge, then $\mu_N$ and $\nu_N$ converge weakly to $\mu$ and $\nu$ as $N \to\infty$, which are absolutely continuous measures with continuous density functions $F(z)$ and $G(z)$, respectively. 
The functions $F(z)$ and $G(z)$ possess continuous partial derivatives of arbitrarily high order, and the functions $F_N(z)$ and $G_N(z)$ and their partial derivatives converge uniformly towards $F(z)$ and $G(z)$ and their partial derivatives as $N \to\infty$. 
\end{lemma}

\begin{lemma}\label{la.3}
Let $f_N(s)$ be a function as in Lemma \ref{la.1}. 
We define measures $\mu_{\sigma,N;\gamma,\delta}$ and $\nu_{\sigma,N;\gamma,\delta}$ as 
\begin{align*}
\mu_{\sigma,N;\gamma,\delta}(A)
&=\frac{1}{\delta-\gamma} \int_{\gamma}^{\delta} 1_A(f_N(\sigma+it)) \,dt, \\
\nu_{\sigma,N;\gamma,\delta}(A)
&=\frac{1}{\delta-\gamma} \int_{\gamma}^{\delta} 1_A(f_N(\sigma+it)) |f'_N(\sigma+it)|^2 \,dt
\end{align*}
for $A \in\mathcal{B}(\mathbb{C})$. 
Suppose that there exist measures $\mu_{\sigma,N}$ and $\nu_{\sigma,N}$ to which $\mu_{\sigma,N;\gamma,\delta}$ and $\mu_{\sigma,N;\gamma,\delta}$ converges weakly as $\delta \to\infty$.  
Suppose further, that $\mu_{\sigma,N}$ and $\nu_{\sigma,N}$ are absolutely continuous with continuous density functions $F_{\sigma,N}(z)$ and $G_{\sigma,N}(z)$, respectively. 

Then the Jensen function $\phi_{f_N}(\sigma)$ is expressed as 
\begin{gather*}
\phi_{f_N}(\sigma)
=\int_{\mathbb{C}} \log|z| F_{\sigma,N}(z) \,|dz|, 
\end{gather*}
and $\phi_{f_N}(\sigma)$ is twice differentiable with the second derivative  
\begin{gather*}
\phi''_{f_N}(\sigma)
=G_{\sigma,N}(0). 
\end{gather*}
\end{lemma}

\begin{proof}
They are direct consequences of the arguments in \cite[Sections 8 and 9]{BorchseniusJessen1948}.
\end{proof}

\begin{proof}[Proof of Theorem \ref{t1.5}]
We take the function $f_N(s)$ in Lemma \ref{la.1} as 
\begin{gather*}
f_N(s)
=\sum_{n=0}^{N} \frac{e^{2\pi i \lambda n}}{(n+\alpha)^s}.
\end{gather*}
Then we have that the limit function $f(s)$ is equal to $L(\lambda,\alpha,s)$, and the assumptions on Lemma \ref{la.1} are satisfied with $\alpha=1/2$,  $\alpha_0>1$, $p=2$ and for any $\beta$, $\beta_0$, and $\gamma_0$. 
Thus if we show the differentiability of the Jensen function $\phi_f(\sigma)$, we have Theorem \ref{t1.5} if we take $\gamma=0$ in formula \eqref{ea.1}. 

Let $\ell(z)=m(z)=z$, $r_n=(n+\alpha)^{-\sigma}$, $\lambda_n=-\log(n+\alpha)$ in Lemma \ref{la.2}. 
Then we have 
\begin{align*}
f_N(s)
&=f_N(e(\gamma_0t), e(\lambda+\gamma_1t), \ldots, e(N\lambda+\gamma_Nt)), \\
f'_N(s)
&=g_N(e(\gamma_0t), e(\lambda+\gamma_1t), \ldots, e(N\lambda+\gamma_Nt)),
\end{align*}
where $e(\theta)=\exp(2\pi i\theta)$ and $\gamma_n=(2\pi)^{-1} \log(n+\alpha)$. 
If $\alpha$ satisfies condition \eqref{c1} of Definition \ref{d1.6}, then $\gamma_0,\ldots,\gamma_N$ are linearly independent over $\mathbb{Q}$. 
Hence we find that for every $z \in \mathbb{C}$, 
\begin{align}\label{eA.2}
&\lim_{T \to\infty} \frac{1}{T} \int_{0}^{T} \psi_z(f_N(\sigma+it)) \,dt\\
&=\lim_{T \to\infty} \frac{1}{T} \int_{0}^{T} 
\psi_z(f_N(e(\gamma_0t), e(\lambda+\gamma_1t), \ldots, e(N\lambda+\gamma_Nt))) \,dt \nonumber\\
&=\int_{T^{N+1}} \psi_z(f_N(\theta_0,\ldots,\theta_N)) \,dm_N(\theta_0,\ldots,\theta_N) \nonumber\\
&=\int_{T^{N+1}} \psi_z(w) \,d\mu_{\sigma,N}(w), \nonumber
\end{align}
where $\mu_{\sigma,N}$ is defined as 
\begin{gather*}
\mu_{\sigma,N}(A)
=\int_{T^{N+1}} 1_A(f_N(t_0,\ldots,t_N)) \,dm_N(t_0,\ldots,t_N)
\end{gather*}
for $A \in\mathcal{B}(\mathbb{C})$. 
Then it is deduced from \eqref{eA.2} that the measure 
\begin{gather*}
\mu_{\sigma,N,T}(A)
=\frac{1}{T} \int_{0}^{T} 1_A(f_N(\sigma+it)) \,dt
\end{gather*}
converges weakly to $\mu_{\sigma,N}$ as $T \to\infty$. 
Similarly, the measure
\begin{gather*}
\nu_{\sigma,N,T}(A)
=\frac{1}{T} \int_{0}^{T} 1_A(f_N(\sigma+it)) |f'_N(\sigma+it)|^2 \,dt
\end{gather*}
converges weakly to
\begin{gather*}
\nu_{\sigma,N}(A)
=\int_{T^{N+1}} 1_A(f_N(t_0,\ldots,t_N)) |g_N(t_0,\ldots,t_N)|^2 \,dm_N(t_0,\ldots,t_N).
\end{gather*}

Then, we use Lemma \ref{la.2}. 
Since $r_n \to0$ as $n \to\infty$, we find that $\mu_{\sigma,N}$ and $\nu_{\sigma,N}$ are absolutely continuous with continuous density functions $F_{\sigma,N}(z)$ and $G_{\sigma,N}(z)$ for large $N$. 
Therefore, by Lemma \ref{la.3}, the Jensen function $\phi_{f_N}(\sigma)$ is expressed as 
\begin{gather*}
\phi_{f_N}(\sigma)
=\int_{\mathbb{C}} \log|z| F_{\sigma,N}(z) \,|dz|, 
\end{gather*}
and $\phi_{f_N}(\sigma)$ is twice differentiable with the second derivative 
\begin{equation}\label{ea.3}
\phi''_{f_N}(\sigma)
=G_{\sigma,N}(0). 
\end{equation} 
Moreover, we see that 
\begin{gather*}
S_0
=\sum_{n=0}^{\infty} \frac{1}{(n+\alpha)^2}, 
\quad
S_1
=\sum_{n=0}^{\infty} \frac{|\log(n+\alpha)|}{(n+\alpha)^2}, 
\quad
S_2
=\sum_{n=0}^{\infty} \frac{\log^2(n+\alpha)}{(n+\alpha)^2}
\end{gather*}
are convergent if $\sigma>1/2$.  
Hence the last part of Lemma \ref{la.2} can be applied, and thus $G_{\sigma,N}(z)$ converges uniformly towards a non-negative continuous function $G_{\sigma}(z)$ as $N \to\infty$. 
Since the Jensen function $\phi_{f_N}(\sigma)$ converges uniformly to $\phi_{f}(\sigma)$ by Lemma \ref{la.1}, we conclude from \eqref{ea.3} that $\phi_{f}(\sigma)$ is twice differentiable with the second derivative  
\begin{gather*}
\phi''_{f}(\sigma)
=G_{\sigma}(0).
\end{gather*}
Therefore Theorem \ref{t1.5} follows from the second part of Lemma \ref{la.1}. 
\end{proof}

\vspace{\baselineskip}

\emph{Acknowledgment.}
The author is deeply grateful to Yoonbok Lee, Kohji Matsumoto, Hidehiko Mishou, Takashi Nakamura, Tomohiro Ooto, J\"{o}rn Steuding, and Masatoshi Suzuki for their helpful comments. 

%\bibliographystyle{amsplain}
%\bibliography{refs}

\providecommand{\bysame}{\leavevmode\hbox to3em{\hrulefill}\thinspace}
\providecommand{\MR}{\relax\ifhmode\unskip\space\fi MR }
\providecommand{\MRhref}[2]{%
  \href{http://www.ams.org/mathscinet-getitem?mr=#1}{#2}
}
\providecommand{\href}[2]{#2}

\end{document}